\documentclass{amsart}

\usepackage{graphicx}
\usepackage{amssymb,mathrsfs}
\usepackage{subfigure}

\usepackage{enumerate}        

\usepackage[all,dvips]{xy}    
\UseComputerModernTips        

\newcommand{\im}{\operatorname{im}}
\newcommand{\dom}{\operatorname{dom}}

\newcommand {\Z}{\mathbb{Z}}
\newcommand {\bng}{B_n\Gamma}

\newcommand {\cng}{\mathcal{C}^n\Gamma}
\newcommand {\ucng}{U\mathcal{C}^n\Gamma}
\newcommand {\ucn}[1]{U\mathcal{C}^n {#1}}
\newcommand {\dng}{\mathcal{D}^n\Gamma}
\newcommand {\udng}{U\mathcal{D}^n\Gamma}
\newcommand {\bn}[1]{B_n #1}
\newcommand {\uc}[2]{U\mathcal{C}^{#1} #2}
\newcommand {\ud}[2]{U\mathcal{D}^{#1} #2}

\newtheorem{theorem}{Theorem}[section]
\newtheorem{lemma}[theorem]{Lemma}

\newtheorem{proposition}[theorem]{Proposition}
\newtheorem{corollary}[theorem]{Corollary}
\newtheorem{conjecture}[theorem]{Conjecture}
\theoremstyle{definition}
\newtheorem{definition}[theorem]{Definition}
\newtheorem{example}[theorem]{Example}

\begin{document}

\title[Cohomology of Tree Braid Groups]{On the Cohomology Rings of tree braid groups}
\author[D.Farley]{Daniel Farley}
      \address{Max Planck Institute for Mathematics\\
               D-53111 Bonn, Germany \\
               http://www.math.uiuc.edu/\~{}farley}
      \email{farley@math.uiuc.edu}
\author[L.Sabalka]{Lucas Sabalka}
      \address{Department of Mathematics\\
               University of Illinois at Urbana-Champaign \\
               Champaign, IL 61820, USA\\
               http://www.math.uiuc.edu/\~{}sabalka}
      \email{sabalka@math.uiuc.edu}

\begin{abstract}
Let $\Gamma$ be a finite connected graph.  The (unlabelled) 
configuration space $\ucng$ of $n$ points on $\Gamma$ is the space of 
$n$-element subsets of $\Gamma$. The $n$-strand braid group of $\Gamma$, 
denoted $\bng$, is the fundamental group of $\ucng$.

We use the methods and results of \cite{FS1} to 
get a partial description of the cohomology rings $H^{\ast}(B_n
T)$, where $T$ is a tree.  Our results are then used to prove that $B_n T$
is a right-angled Artin group if and only if $T$ is linear or $n<4$.  
This gives a large number of counterexamples to Ghrist's conjecture that
braid groups of planar graphs are right-angled Artin groups.
\end{abstract}

\keywords{tree braid groups, configuration spaces, right-angled Artin groups}

\subjclass[2000]{Primary 20F65, 20F36; Secondary 57M15, 55R80}

\maketitle

\section{Introduction}\label{sec:intro}

If $\Gamma$ is a finite connected graph and $n$ is a natural number, 
then the \emph{unlabelled configuration space} of $n$ points on 
$\Gamma$, denoted $\ucng$, is the space of $n$-element subsets of 
$\Gamma$, endowed with the Hausdorff topology.  The \emph{labelled 
configuration space} $\cng$ is the space of $n$-tuples of distinct 
elements in $\Gamma$.  The \emph{$n$-strand braid group of $\Gamma$}, 
denoted $\bng$, is the fundamental group of $\ucng$; the 
\emph{$n$-strand pure braid group of $\Gamma$}, $P \bng$, is the 
fundamental group of $\cng$.

Various properties of graph braid groups have been established by other
authors.  Ghrist showed in \cite{Ghrist} that the spaces $\cng$ are $K(
PB_n \Gamma, 1)$s, and that a $K ( PB_n \Gamma, 1)$ is homotopy equivalent
to a complex of dimension at most $k$, where $k$ is the number of vertices
in $\Gamma$ of degree at least $3$.  He also made the following
conjecture:

\begin{conjecture} \label{conj:Ghrist}
\cite{Abrams,Ghrist}  The (pure) braid group of any planar graph is a 
right-angled Artin group.
\end{conjecture}   

Abrams \cite{Abrams} (for all $n$) and Hu \cite{Hu} (for the case $n = 
2$) introduced a discretized configuration space $\dng$, and showed that 
$\cng$ and $\dng$ are homotopy equivalent under appropriate hypotheses 
(which are easy to satisfy). Abrams went on to prove that the universal 
cover of the space $\dng$ is a CAT(0) cubical complex.  This implies, in 
particular, that graph braid groups have solvable word and conjugacy 
problems \cite{BridsonHaeflinger}.  Abrams also showed that $PB_2 ( K_5 
)$ and $PB_2 ( K_{3,3} )$ are the fundamental groups of closed surfaces, 
and thus aren't right angled Artin groups.  This is the reason for the 
word ``planar" in Conjecture \ref{conj:Ghrist}.  Crisp and Wiest 
\cite{CrispWiest} have shown that all graph braid groups embed in 
right-angled Artin groups, which implies that graph braid groups are 
linear, bi-orderable, residually finite, and residually nilpotent. 
Connolly and Doig \cite{ConnollyDoig} showed that the braid group of any 
linear tree is a right-angled Artin group.  (A tree $T$ is \emph{linear} 
if there is an embedded arc which passes through every vertex in $T$ of 
degree at least $3$.)

This paper continues a project begun in \cite{FS1}.  In \cite{FS1}, we
used a discrete version of Morse theory (due to Forman \cite{Forman}) to
simplify the configuration spaces $\ucng$ within their homotopy types.  
Our immediate goal was to settle Conjecture \ref{conj:Ghrist}.  We were
able to compute presentations $\mathcal{P}( B_n T )$ for all braid groups
$B_n T$, where $T$ is a tree; that is, for all \emph{tree braid groups}
(\cite{FS1}, Theorem 5.3).  The generators of $\mathcal{P}( B_n T )$ are
in one-to-one correspondence with critical $1$-cells of $\ucng$ and
relators correspond to critical $2$-cells.  Here ``critical" is used in
the sense of Forman's discrete Morse theory.  In \cite{Farley2} it was
shown that $H_{i} ( \ucn{T})$ (equivalently, $H_{i} ( B_n T)$, since
$\ucng$ is aspherical for any graph $\Gamma$ \cite{Abrams}, \cite{Ghrist})
is a free abelian group of rank equal to the number of critical $i$-cells
in $\ucn{T}$.  It follows from this that $\mathcal{P}( B_n T)$ has the
minimum possible number of generators and relators.  We were unable to
produce counterexamples to Conjecture \ref{conj:Ghrist}, although the form
of the relators in $\mathcal{P}( B_n T)$ made a negative answer seem
likely for most trees and most natural numbers $n$.

Here we get nearly complete information about the mod $2$ cohomology 
rings of tree braid groups.  Our results allow us to prove that most 
tree braid groups are not right-angled Artin groups (see Theorem 
\ref{thm:Ghrist}). Thus we produce a large number of counterexamples to 
the version of Conjecture \ref{conj:Ghrist} in which the word ``pure" is 
omitted.  (It is worth noting here that Abrams and Ghrist made the 
conjecture only for pure braid groups.  In this paper, we refer to 
either version of Conjecture \ref{conj:Ghrist} as ``Ghrist's 
conjecture".  We believe that the analogue of Theorem \ref{thm:Ghrist} 
will be true for pure braid groups.)

The argument is as follows.  We first compute the cohomology ring of $B_4
T_{min}$, where $T_{min}$ is the minimal nonlinear tree.  Our calculation
shows that $B_4 T_{min}$ is not a right-angled Artin group, since
$H^{\ast}( B_4 T_{min} ; \mathbb{Z}/2\mathbb{Z} )$ is not the exterior
face ring of a flag complex (see Section \ref{sec:Counterex}).  If $T$ is
any nonlinear tree and $n \geq 4$, we embed $\uc{4}{T_{min}}$ into
$\uc{n}{T}$.  By analyzing the kernel of the map on cohomology, we can
conclude that $B_n T$ is also not a right-angled Artin group, since its
cohomology ring also fails to be the exterior face ring of a flag complex.

Finally, we note that our description of the mod $2$ cohomology rings of
tree braid groups is likely to have other applications.  For instance,
Michael Farber \cite{Farber1}, \cite{Farber2} has defined an invariant
$TC(X)$ of a topological space $X$, called the \emph{topological
complexity} of $X$, which is an integer measuring the complexity of
motion-planning problems of systems having $X$ as their configuration
space.  Farber establishes cohomological lower bounds for $TC(X)$ in
\cite{Farber1}.

This paper is organized as follows.  In Section \ref{sec:background} we
give a brief description of discrete Morse theory and its applications to
computing homology.  In Section \ref{sec:DGVF} we describe 
Morse matchings on the spaces $\ucng$.  In Section
\ref{sec:mod2ring}, we give a partial description of the mod $2$
cohomology ring of any tree braid group.  In Section \ref{sec:Counterex},
we use the results of Section \ref{sec:mod2ring} and a cohomological
argument to determine which tree braid groups are right-angled Artin
groups.

We would like to thank Ilya Kapovich and Robert Ghrist for participating 
in discussions related to this work.  We thank Aaron Abrams and Carl 
Mautner for telling us of some of Mautner's counterexamples to 
Conjecture \ref{conj:Ghrist} (in both the pure and regular cases) 
\cite{Mautner}, and for suggesting that cohomological methods might be 
used to produce counterexamples.

\section{Background on Discrete Morse Theory}\label{sec:background}

\subsection{Basic Definitions} \label{sec:morsetheory} 
In this subsection, we collect some basic definitions from \cite{FS1} (see
also \cite{Brown} and \cite{Forman}, which were the original sources for
these ideas).

Let $X$ be a finite regular CW complex. Let $K$ denote the set of open
cells of $X$.  Let $K_{p}$ be the set of open $p$-cells of $X$.  For open
cells $\sigma$ and $\tau$ in $X$, we write $\sigma < \tau$ if $\sigma \neq
\tau$ and $\sigma \subseteq \overline{\tau}$, where $\overline{\tau}$ is
the closure of $\tau$, and $\sigma \leq \tau$ if $\sigma < \tau$ or
$\sigma = \tau$.

A \emph{partial function} from a set $A$ to a set $B$ is a function
defined on a subset of $A$, and having $B$ as its target.  A
\emph{discrete vector field} $W$ on $X$ is a sequence of partial functions
$W_{i}: K_{i} \rightarrow K_{i+1}$ such that:

\begin{enumerate}
\item Each $W_{i}$ is injective;
\item if $W_{i} (\sigma ) = \tau$, then $\sigma < \tau$;
\item $\im \left( W_{i} \right) \cap \dom \left( W_{i+1} \right) = 
\emptyset$, where $\im$ denotes image and $\dom$ denotes domain.
\end{enumerate}

Let $W$ be a discrete vector field on $X$.  A \emph{$W$-path of 
dimension $p$} is a sequence of $p$-cells $\sigma_{0}, \sigma_{1}, 
\ldots, \sigma_{r}$ such that if $W( \sigma_{i} )$ is undefined, then 
$\sigma_{i+1} = \sigma_{i}$; otherwise $\sigma_{i+1} \neq \sigma_{i}$ 
and $\sigma_{i+1} < W( \sigma_{i})$.  The $W$-path is \emph{closed} if 
$\sigma_{r} = \sigma_{0}$, and \emph{non-stationary} if $\sigma_{1} \neq 
\sigma_{0}$.  A discrete vector field $W$ is a \emph{Morse matching} if 
$W$ has no non-stationary closed paths.  

If $W$ is a Morse matching, then a cell $\sigma \in K$ is 
\emph{redundant} if it is in the domain of $W$, \emph{collapsible} if it 
is in the image of $W$, and \emph{critical} otherwise.  Note that any 
two of these categories are mutually exclusive by condition (3) in the 
definition of discrete vector field.

The ideas ``discrete Morse function" and ``Morse matching" are largely 
equivalent, in a sense that is made precise in \cite{Forman}, pg. 131.  
In practice, we will always use Morse matchings instead of discrete 
Morse functions in this paper (as we also did in \cite{FS1}). A Morse 
matching is sometimes referred to as a ``discrete gradient vector 
field'' in the literature; in particular, this is the case in 
\cite{FS1}.

\subsection{Discrete Morse Theory and Homology}
The discrete Morse theory sketched in Subsection \ref{sec:morsetheory} can
be used to compute homology groups.  We include only a brief account,
without proofs.  More extended expositions can be found in \cite{Forman}
and \cite{Farley1}.

Fix an oriented finite regular CW complex $X$.  Let $C_{\ast}(X)$ be the 
cellular chain complex of $X$.  Each chain group $C_{n}(X)$ has a 
distinguished basis consisting of positively oriented $n$-cells, denoted 
$B_{n}(X)$.  Let $W$ be a Morse matching, and define a 
map $\widehat{W}_{n} : C_{n}(X) \rightarrow C_{n+1}(X)$ as follows:

\begin{eqnarray*}
\widehat{W}_{n} (c) & = & \pm W(c)  \hbox{ if $c$ is redundant;} \\
\widehat{W}_{n} (c) & = & 0  \hbox{ otherwise.}
\end{eqnarray*}

Here the sign is chosen so that the oriented cell $c$ occurs with the
coefficient $-1$ in $\partial \widehat{W}_{n} (c)$ if $c$ is redundant.  
Extend linearly to a map $\widehat{W}_{n}: C_{n}(X) \rightarrow
C_{n+1}(X)$.  Define a chain map $f_{\widehat{W}}: C_{\ast}(X) \rightarrow
C_{\ast}(X)$, called the \emph{discrete flow} associated to $\widehat{W}$,
by setting $f_{\widehat{W}} = 1 + \partial \widehat{W} + \widehat{W}
\partial$.  We usually omit the subscript and simply write $f$.

The discrete flow $f$ has the following properties:

\begin{lemma} \label{lem:homology}  ~\newline 
\begin{enumerate}
\item (\cite{Farley1}, \cite{Forman}) For any finite chain $c \in
C_{\ast}(X)$, there is some $m \in \mathbb{N}$ such that $f^{m}(c) =
f^{m+1}(c) = \ldots$.  It follows that there is a well-defined chain map
$f^{\infty}: C_{\ast}(X) \rightarrow C_{\ast}(X)$.

\item (\cite{Farley1}; cf. \cite{Forman}) If $c$ is any cycle in
$C_{\ast}(X)$, then there is a unique $f$-invariant cycle that is
homologous to $c$, namely $f^{\infty}(c)$. Moreover, $f^{\infty}(c)$ is a
linear combination of oriented critical cells and collapsible cells (i.e.,
any redundant cell appears with a coefficient of $0$).

\item (\cite{Farley1}) If $c$ is a collapsible cell, then $f^{\infty}(c) =
0$.  If $c$ is critical, then $f^{\infty}(c) = c + (collapsible~cells)$.  
As a result, an $f$-invariant chain is determined by its critical cells,
i.e., if $c$ is an $f$-invariant chain and $c = c_{crit} + c_{coll}$,
where $c_{crit}$ is a linear combination of critical cells and $c_{coll}$
is a linear combination of collapsible cells, then
  $$ c = f^{\infty}(c) = f^{\infty}(c_{crit}).$$
\end{enumerate}
\end{lemma}
                                                     
Properties (1)-(3) show that if a finite regular CW complex $X$ is 
endowed with a Morse matching $W$, then the homology groups of $X$ are 
largely determined by the critical cells of $X$.  We now make this 
statement more precise.  Fix a Morse matching $W$.  For 
$i \geq 0$, let $M_{i}(X)$ denote the free abelian group on the set of 
positively oriented critical $i$-cells.  Give the collection of abelian 
groups $M_{i}(X) \, \, (i \geq 0)$ the structure of a chain complex, 
called the \emph{Morse complex}, by identifying $M_{i}(X)$ with 
$C_{i}(X)$ via the map $f^{\infty}$.  The boundary map 
$\tilde{\partial}$ in the Morse complex is defined by
  $$ \tilde{\partial}(c) = \Pi \partial f^{\infty}(c)  \quad 
  \left(  c \in M_{i}(X) \right),$$
where $\Pi$ denotes projection onto the factor of $C_{i-1}(X)$ spanned by
the critical $i-1$ cells.

We have the following theorem:

\begin{theorem} \emph{(}\cite{Farley1}, \cite{Forman}\emph{)}
The Morse complex $\left( M_{n}(X), \tilde{\partial}_{n} \right)$ and the
cellular chain complex $\left( C_{n}(X), \partial_{n} \right)$ have
isomorphic homology groups, by an isomorphism which sends a cycle $c$ from
the Morse complex to $f^{\infty}(c)$.
\qed
\end{theorem}

\section{A Morse Matching on the Discretized Configuration
Space $UD^{n} \Gamma$} \label{sec:DGVF}

\subsection{Definitions and an Example} \label{sec:definitions}
Throughout this paper, all graphs are assumed to be finite and 
connected.

Let $\Gamma$ be a graph, and fix a natural number $n$. The \emph{labelled
configuration space} of $\Gamma$ on $n$ points is the space
  $$\left(\prod^{n} \Gamma \right) - \Delta,$$
where $\Delta$ is the set of all points $(x_1, \dots, x_n) \in \prod^n
\Gamma$ such that $x_i = x_j$ for some $i \neq j$. The \emph{unlabelled
configuration space} of $\Gamma$ on $n$ points is the quotient of the
labelled configuration space by the action of the symmetric group $S_n$,
where the action permutes the factors. The \emph{braid group} of $\Gamma$
on $n$ strands, denoted $B_n\Gamma$, is the fundamental group of the
unlabelled configuration space of $\Gamma$ on $n$ strands. The \emph{pure
braid group}, denoted $PB_n\Gamma$, is the fundamental group of the
labelled configuration space.

The set of vertices of $\Gamma$ will be denoted by $V(\Gamma)$, and the
degree of a vertex $v \in V(\Gamma)$ is denoted $d(v)$.  If a vertex $v$
is such that $d(v) \geq 3$, $v$ is called \emph{essential}.

Let $\Delta'$ denote the union of those open cells of $\prod^n \Gamma$
whose closures intersect $\Delta$. Let $\dng$ denote the
space $\prod^n \Gamma - \Delta'$. Note that $\dng$ inherits a CW complex
structure from the Cartesian product, and that a cell in $\dng$ has the
form $c_1 \times \dots \times c_n$ such that each $c_i$ is either a vertex
or the interior of an edge, and the closures of the $c_i$ are mutually
disjoint. Let $\udng$ denote the quotient of $\dng$ by the action of the
symmetric group $S_n$ which permutes the coordinates. Thus, an open cell in
$\udng$ has the form $\{c_1, \dots, c_n\}$ such that each $c_i$ is either
a vertex or the interior of an edge and the closures are mutually
disjoint.  The set notation is used to indicate that order does not
matter.

Under most circumstances, the labelled (respectively, unlabelled)
configuration space of $\Gamma$ is homotopy equivalent to $\dng$
(respectively, $\udng$).  Specifically:

\begin{theorem} \cite{Abrams} \label{thm:Abrams}
For any $n>1$ and any graph $\Gamma$ with at least $n$ vertices, the
labelled (unlabelled) configuration space of $n$ points on $\Gamma$ strong
deformation retracts onto $\dng$ ($\udng$) if
\begin{enumerate}
\item each path between distinct vertices of degree not equal to $2$
passes through at least $n-1$ edges; and

\item each path from a vertex to itself which is not null-homotopic in
$\Gamma$ passes through at least $n+1$ edges.
\end{enumerate}
\end{theorem}

A graph $\Gamma$ satisfying the conditions of this theorem for a given $n$
is called \emph{sufficiently subdivided} for this $n$.  It is clear that
every graph is homeomorphic to a sufficiently subdivided graph, no matter
what $n$ may be.

Throughout the rest of the paper, we work exclusively with the space
$\udng$ where $\Gamma$ is sufficiently subdivided for $n$.  Also from now
on, ``edge'' and ``cell'' will refer to closed objects.

Choose a maximal tree $T$ in $\Gamma$. Edges outside of $T$ are called
\emph{deleted edges}.  Pick a vertex $\ast$ of valence $1$ in $T$ to be
the root of $T$. Choose an embedding of the tree $T$ into the plane.  We
define an order on the vertices of $T$ (and, thus, on vertices of
$\Gamma$) as follows. Begin at the basepoint $\ast$ and walk along the
tree, following the leftmost branch at any given intersection, and
consecutively number the vertices in the order in which they are first
encountered.  (When you reach a vertex of degree one, turn around.)  The
vertex adjacent to $\ast$ is assigned the number $1$.  Note that this
numbering depends only on the choice of $\ast$ and the embedding of the
tree.  Let $\iota(e)$ and $\tau(e)$ denote the endpoints of a given edge
$e$ of $\Gamma$.  Without loss of generality, we orient each edge to go
from $\iota(e)$ to $\tau(e)$, and so that $\iota(e) > \tau(e)$.  (Thus, if
$e \subseteq T$ the geodesic segment $[ \iota(e), \ast ]$ in $T$ must pass
through $\tau(e)$.)

We use the order on the vertices to define a Morse matching $W$ on 
$UD^{n} \Gamma$.  We begin with some definitions which will help to 
classify cells of $UD^{n} \Gamma$ as critical, collapsible, or 
redundant.

Let $c = \{ c_1, \ldots, c_{n-1}, v \}$ be a cell in $\ud{n}{\Gamma}$ 
containing a vertex $v$.  If $v = \ast$, then $v$ is \emph{blocked} in 
$c$; otherwise, let $e$ be the unique edge in $T$ such that $\iota(e) = 
v$.  If $e \cap c_i \neq \emptyset$ for some $i \in \{1, \dots, n-1\}$, 
we also say $v$ is \emph{blocked} in $c$; otherwise, $v$ is unblocked.  
Equivalently, $v$ is unblocked in $c$ if and only if $\{ c_1, \ldots, 
c_{n-1}, e \}$ is also a cell in $\ud{n}{\Gamma}$.  If $c = \{ c_1, 
\ldots, c_{n-1}, e \}$, the edge $e$ is \emph{disrespectful} 
in $c$ if \begin{enumerate}
  \item there is a vertex $v$ in $c$ such that 
  \begin{enumerate}
    \item $v$ is adjacent to $\tau(e)$, and
    \item $\tau(e) < v < \iota(e)$, or
  \end{enumerate}    
  \item $e$ is a deleted edge.
\end{enumerate}
Otherwise, the edge $e$ is \emph{respectful} in $c$.  Conceptually, 
think of the edge $e$ as representing a strand in $c$ moving from 
$\iota(e)$ to $\tau(e)$.  Then $e$ is disrespectful in $c$ if that 
strand is moving out of turn by not respecting the order on vertices of 
$T$.  In the paper \cite{FS1}, disrespectful was referred to by ``non 
order-respecting''.

It will occasionally be useful to have another definition.  If $v$ is a 
vertex in the tree $T$, we say that two vertices $v_1$ and $v_2$ lie in 
the same \emph{direction} from $v$ if the geodesics $[ v , v_1 ] , [v , 
v_2 ] \subseteq T$ start with the same edge.  Thus, there 
are $deg(v)$ directions from a vertex of degree $deg(v)$ in $T$.  We 
number these directions $0, 1 , 2, \ldots, deg(v)-1$, beginning with the 
direction represented by $[v, \ast]$, numbered $0$, and proceeding in 
clockwise order.  We will sometimes write $g( v_1 , v_2 )$ (where $v_1 
\neq v_2$) to refer to the direction from $v_1$ to $v_2$.

Suppose that we are given a cell $c = \{ c_1, \ldots, c_n \}$ in $UD^{n}
\Gamma$.  Assign each cell in $c$ a number as follows.  A vertex of $c$ is
given the number from the above traversal of $T$.  An edge $e$ of $c$ is
given the number for $\iota(e)$.  Arrange the cells of $c$ in a sequence
$\mathcal{S}$, from the least- to the greatest-numbered.  The following
definition of a Morse matching $W$ is equivalent to the
definition of $W$ from \cite{FS1}, by Theorem 3.6 of the same paper.

\begin{definition} \label{def:critical}
We define a Morse matching $W$ on $\ud{n}{\Gamma}$ 
as follows:
\begin{enumerate}
\item If an unblocked vertex occurs in $\mathcal{S}$ before all of the
respectful edges in $c$ (if any), then $W(c)$ is obtained from $c$
by replacing the minimal unblocked vertex $v \in c$ with $e(v)$, where
$e(v)$ is the unique edge in $T$ satisfying $\iota (e(v)) = v$.  In
particular, $c$ is redundant.

\item If a respectful edge occurs before any unblocked vertex, then
$c \in \im W$, i.e., $c$ is collapsible.  The cell $W^{-1}(c)$ is
obtained from $c$ by replacing the minimal respectful edge $e$ with
$\iota(e)$.

\item If there are neither unblocked vertices nor respectful edges
in $c$, then $c$ is critical.
\end{enumerate}
\end{definition}

\begin{example}
Figure \ref{fig:example} depicts three different cells of $UD^{4} T_{min}$
for the given tree $T_{min}$.  In each case, the vertices and edges of the
given cell are numbered from least to greatest, in the sense mentioned
above.  (The numbering of these cells differs from the above-described
order, but this doesn't matter since the ordering remains the same.  For
instance, the vertices and edges in the cell pictured in Figure
\ref{fig:example} (a) should be numbered $10$, $14$, $16$, and $19$,
instead of (respectively) $1$, $2$, $3$, $4$.)
\begin{figure}[!h]
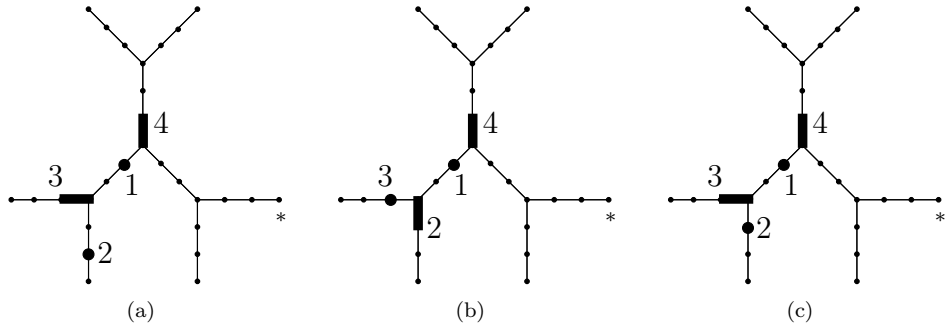

\centering

\subfigure[]{
  \input{example1.pstex_t}
}\hfill
\subfigure[]{
  \input{example2.pstex_t}
}\hfill
\subfigure[]{
  \input{example3.pstex_t}
}

\caption{Three different cells of $\ud{n}{T}$.}
\label{fig:example}
\end{figure}

The vertex numbered $1$ in (a) is blocked.  The vertex numbered $2$ is
unblocked, so the cell in (a) is redundant.  Note that edge $3$ is
respectful, and edge $4$ is disrespectful.  We get $W(c_1)$
by replacing vertex $2$ with the unique edge in $T$ having vertex $2$ as
its initial vertex.  In terms of the usual ordering, this is the edge $[
v_{14}, v_{13} ]$.

Let $c_2$ denote the cell depicted in (b).  The vertex numbered $1$ in
$c_2$ is blocked.  The edge numbered $2$ is respectful, so $c_2$ is
collapsible.  Note that vertex $3$ is blocked and edge $4$ is 
disrespectful.  The description of $W^{-1}$ above implies that
$W^{-1}(c_2)$ is obtained from $c_2$ by replacing edge $2$ with its
initial vertex.

The cell depicted in (c) is critical since vertices $1$ and $2$ are both
blocked, and the edges $3$ and $4$ are disrespectful.
\end{example}

\section{The Mod $2$ Cohomology Ring of $\ud{n}{T}$}\label{sec:mod2ring}

In this section, we give a partial description of the mod $2$ cohomology
ring of $\ud{n}{T}$, where $n$ is an arbitrary natural number and $T$ is
an arbitrary tree. The method is to map $\ud{n}{T}$ to a new complex
$\widehat{\ud{n}{T}}$, which is a subcomplex of a high-dimensional torus.  
The induced map $q^{\ast} : H^{\ast} (\widehat{\ud{n}{T}}) \rightarrow
H^{\ast}( \ud{n}{T} )$ turns out to be surjective.  This gives us an easy
way to compute the cup product:  we take two cohomology classes in
$H^{\ast} ( \ud{n}{T})$, look at their preimages under the map $q^{\ast}$,
cup these preimages using known facts about the cohomology rings of
subcomplexes of torii, and then push the product back over into $H^{\ast}(
\ud{n}{T})$.

The complex $\widehat{\ud{n}{T}}$ can be described very simply:  it is the
result of identifying the opposite sides of all of the cubes in
$\ud{n}{T}$.  Thus $\widehat{\ud{n}{T}}$ consists of a union of
(potentially singular) torii, one for each cell of $\ud{n}{T}$.

It is by no means clear, however, and false in general, that identifying
all opposite faces in a CAT(0) cubical complex will result in a subcomplex
of a torus.  For this reason, we give a very careful (and somewhat
abstract) proof that $\widehat{\ud{n}{T}}$ has the properties we want.

\subsection{An Equivalence Relation on the Cells of $\ud{n}{T}$}

Let $E(c)$ denote the set of edges of the $i$-cell $c$.  We abuse the 
notation by also letting $E(c)$ denote the subset $\bigcup_{e \in E(c)} 
e$ of $T$.

Let $c$ and $c'$ be $i$-cells $(0 \leq i \leq n)$ of $\ud{n}{T}$.  Say 
$c = \{ e_1 , \ldots , e_i , v_{i+1}, \ldots, v_n \}$ and $c' = \{ 
e'_{1} , \ldots , e'_{i}, v'_{i+1}, \ldots, v'_n \}$, where $E(c) = 
\{e_1, \dots, e_i\}$ and $E(c') = \{e'_1, \dots, e'_i\}$.  Write $c \sim
> c'$ if

\begin{enumerate}
\item $E(c) = E(c')$, and
\item for any connected component $C$ of $T - E(c)$,
$$ \left| C \cap \left\{ v_{i+1}, \ldots, v_{n} \right\} \right| = 
   \left| C \cap \left\{ v'_{i+1}, \ldots, v'_{n} \right\} \right|.$$
\end{enumerate}
Let $K$ be the set of open cells in $\ud{n}{T}$, as in Section
\ref{sec:morsetheory}.  It is rather clear that $\sim$ is an equivalence
relation on $K$. Let $[c]$ denote the equivalence class of a cell $c$.

We define a partial order $\leq$ on the equivalence classes based on the
partial order on cells, writing $[c] \leq [c_1]$ if there exist
representatives $\hat{c} \in [c]$, $\hat{c}_1 \in [c_1]$ such that
$\hat{c} \leq \hat{c}_1$. It is slightly nontrivial to verify that $\leq$
is transitive on the equivalence classes. Suppose that $[c_1] \leq [c_2]$
and $[c_2] \leq [c_3]$.  There are representatives $\hat{c}_1 \in [c_1],
\hat{c}_2, \tilde{c}_2 \in [c_2]$, and $\tilde{c}_3 \in [c_3]$ such that
$\hat{c}_1 \leq \hat{c}_2$ and $\tilde{c}_2 \leq \tilde{c}_3$.  Since
$\hat{c}_2 \sim \tilde{c}_2$, it is possible, by moving the vertices of
$\hat{c}_2$ one at a time along edges of $T$, and leaving all edges fixed,
to arrive at $\tilde{c}_2$.  Every vertex of $\hat{c}_2$ is also a vertex
of $\hat{c}_1$; if a vertex of $\hat{c}_1$ is also in $\hat{c}_2$, then
move it as above.  Call the result of doing these moves $\tilde{c}_1$.  
Then $\tilde{c}_1 \sim \hat{c}_1$ and it is clear that $\tilde{c}_1 \leq 
\tilde{c}_2 \leq \tilde{c}_3$, so transitivity of $\leq$ follows.

We state the following lemma in terms of a tree $T$, but we note that 
with the appropriate definitions the statements in parts (1), (4), and 
(5) may be generalized to arbitrary graphs.

\begin{lemma}[Properties of $\sim$ and $\leq$] \label{lem:big}
Let $T$ be a finite connected tree.
\begin{enumerate}
\item If $c$ and $c'$ are $i$-cells of $\ud{n}{T}$ with $E(c) = 
E(c')$ and the equivalence classes $[c],[c']$ have a common upper 
bound $[\tilde{c}]$ with respect to $\leq$, then $[c] = [c']$.

\item Let $c_1, \ldots, c_j$ be $1$-cells in $\ud{n}{T}$ from distinct 
equivalence classes. If $\{ [c_1], \ldots, [c_j] \}$ has an upper bound 
with respect to $\leq$, then $\{ [c_1], \ldots, [c_j] \}$ has a least 
upper bound with respect to $\leq$.  Furthermore, if $e_1, \ldots, e_j$ 
are edges of $T$ satisfying $e_i \in c_i$ for $1 \leq i \leq j$, then 
$e_1, \ldots, e_j$ are pairwise disjoint.

\item If $\tilde{c}$ is a $j$-cell in $\ud{n}{T}$, then there is a 
unique collection $\{ [c_1], \ldots, [c_j] \}$ of equivalence classes of 
$1$-cells such that $[\tilde{c}]$ is the least upper bound of $\{[c_1], 
\ldots, [c_j] \}$ with respect to $\leq$.

\item If $c$ is a critical cell in $\ud{n}{T}$ and $[c'] \leq [c]$,
then $c' \sim \hat{c}$ for some critical cell $\hat{c}$.

\item If $c$ is a critical cell in $\ud{n}{T}$ and $c' \sim c$, then
$c' = c$ or $c'$ is redundant.
\end{enumerate}
\end{lemma}

\begin{proof}
(1)  Let $[\tilde{c}]$ be an equivalence class of $j$-dimensional cells in
$\ud{n}{T}$, where $\tilde{c} = \{ e_1, \ldots, e_j, v_{j+1}, \ldots,
v_n\}$.  It is enough to show that, for any $c$ with $[c] \leq
[\tilde{c}]$, the equivalence class $[c]$ is uniquely determined by the
collection $\{ e_{i_1}, \ldots, e_{i_{k}} \}$ of all edges common to both
$\tilde{c}$ and $c$.

An arbitrary face $c$ of $\tilde{c}$ is determined by selecting a subset of
the edges $\{ e_1, \ldots, e_j \}$ and replacing each edge of this subset
by either its initial or its terminal vertex.  By an argument similar to
that establishing
the transitivity of $\leq$, the equivalence class of $c$ depends neither
on the representative chosen from $[\tilde{c}]$ nor on the choice involved
in replacing an edge with one of its endpoints.  Thus, given $\tilde{c}$
with $[c] \leq [\tilde{c}]$, $[c]$ is uniquely determined by the edges $c$
has in common with $\tilde{c}$ - i.e. $[c]$ is uniquely determined by
$E(c)$.  This proves part (1).

(2)
Suppose that $[c]$ is an upper bound for $\{ [c_1], \ldots, [c_j] \}$
where the $[c_1], \ldots, [c_j]$ are all distinct.  For $i \in \{ 1,
\ldots, j \}$, let $e_i$ be the unique edge in $c_i$.  Certainly, $\{ e_1,
\ldots, e_j \} \subseteq c$.  Let $c'$ be a cell given by replacing any
extra edges $e \in c - \{ e_1, \ldots, e_j \}$ with either $\iota(e)$ or
$\tau(e)$.  Then $[c_i] \leq [c']$ for $i \in \{ 1, \ldots, j \}$,
so $[c']$ is an upper bound for $\{ [c_1], \ldots, [c_j] \}$.  It follows
that given any upper bound $[c]$ for $\{ [c_1], \ldots, [c_j] \}$
there exists another upper bound $[c']$ such 
that $[c'] \leq [c]$ and $E(c') = \{ e_1, \ldots, e_j \}$.

To prove the first claim of part (2), it remains to be shown that there is
only one upper bound $[c']$ such that $E(c') = \{ e_1 , \ldots , e_j \}$.  
Suppose that $[c']$ and $[c'']$ are
both such upper bounds for $\{ [c_1], \ldots, [c_j] \}$.  Thus $E(c') = E(c'')
= \{e_1, \ldots, e_j\}$.  Fix an integer $i$, $1 \leq i \leq j$.  Let $C$
be a connected component of $T-e_i$, and assign to $C$ the integer
$f_{e_i}(C)$ defined by:
\begin{eqnarray*}
  \hspace*{1cm} f_{e_i}(C) & = & \left( \#\{\hbox{vertices or edges of }c' 
  \hbox{ contained in } C\}\right) \\ 
  & & - \left( \#\{\hbox{vertices or edges of }c'' \hbox{ contained in } C\}\right).
\end{eqnarray*}  

It must be that $f_{e_i}(C) = 0$ for every $i$ and every component $C$ of
$T-e_i$.  For if not, then for some $i$ there exist distinct equivalence
classes $[c_0'], [c_0'']$ of $1$-cells with $E(c_0') = E(c_0'') = \{e_i\}$
such that $[c_0'] \leq [c']$ and $[c_0''] \leq [c'']$.  But by part (1),
we must have that $[c_i] = [c_0']$ since $[c']$ is a common upper bound for
$[c_i]$ and $[c_0']$.  Similarly, $[c_i] = [c_0'']$.  This means that $[c_0']
= [c_0'']$, which is a contradiction, so indeed $f_{e_i}(C) = 0$ for any
$i$ and any $C$.

If $[c'] \neq [c'']$, then $c'$ and $c''$ have a different number of 
cells in some connected component of $T - E(c')$.  Fix a connected 
component $C$ of $T - E(c')$.  Let $e_{i_1}, \ldots, e_{i_l}$ be the 
edges of the collection $\{e_1, \ldots, e_j \}$ that have exactly one 
endpoint in the component $C$.  We analyze the connected components of 
$T - \bigcup_{k = 1}^l e_{i_k}$.  These connected components are either 
$C$ itself, or contain a single endpoint from exactly one of the edges 
$e_{i_1}, \ldots, e_{i_l}$. (Note: here we have just used the fact that 
$T$ is a tree for the first time.)  If a connected component $C'$ of $T 
- \bigcup_{k = 1}^l e_{i_k}$ is not $C$ itself, and contains an endpoint 
of $e_{i_m}$, say, then it contains equal numbers of cells from $c'$ and 
$c''$ by the claim, since $C'$ is in fact a connected component of $T - 
e_{i_m}$.  (This again uses the fact that $T$ is a tree.) It follows by 
process of elimination that $C$ contains equal numbers of cells, 
necessarily vertices, from both $c'$ and $c''$.  Since $C$ was an 
arbitrary connected component of $T - E(c')$, it must be that $[c'] = 
[c'']$.

We now prove the second claim.  Suppose that $[c_1], \ldots,
[c_j]$ are distinct equivalence classes of $1$-cells having a common upper
bound $[c]$.  Part (1) shows that $e_{i_1} \neq e_{i_2}$ for two distinct
equivalence classes $[c_{i_1}], [c_{i_2}] \in \{ [c_1], \ldots, [c_j] \}$.  
If $e_{i_1} \in c_{i_1}$, $e_{i_2} \in c_{i_2}$, $e_{i_1} \cap e_{i_2}
\neq \emptyset$, then $e_{i_1}, e_{i_2} \in c$, which is impossible since
$c$ is a cell of $\ud{n}{T}$ and $e_{i_1} \cap e_{i_2} \neq \emptyset$.  

(3)  Let $S = \{ [c] \mid \dim c = 1 \hbox{ and }
[c] \leq [\tilde{c}] \}$.  Part (1) implies that an element of $S$ is
uniquely determined by a choice of edge from $\tilde{c}$.  Thus $|S| = j$.  
The fact that $[\tilde{c}]$ is the least upper bound of $S$ follows from 
the description of the least upper bound in (2).

(4) If $c$ is a critical cell of $\ud{n}{\Gamma}$ and $[c_1] < [c]$, 
then there is some representative of the equivalence class $[c_1]$ - say 
$c_1$ - such that $c_1 < c$.  For, let $c_1' \in [c_1]$ and $c' \in [c]$ 
be such that $c_1' < c'$.  Then $c_1$ (respectively, $c$) is the result 
of moving all vertices in $c_1'$ (respectively, $c'$) toward $\ast$ 
until they are blocked.  Thus each edge in $c_1$ is an edge in $c$, and 
each vertex in $c$ is a vertex in $c_1$.  It follows that no edges in 
$c_1$ are respectful.  Now repeatedly move each unblocked vertex 
of $c_1$ toward $\ast$ until it is blocked.  This operation clearly 
preserves $\sim$, and the resulting cell is critical, having no 
unblocked vertices and no respectful edges.

(5) Suppose $c$ is critical and $c_1 \sim c$.  Suppose first that $c_1$ 
has no respectful edges.  If $c_1$ has unblocked vertices, then it 
follows that $c_1$ is redundant.  If $c_1$ has no unblocked vertices, 
then it is critical by Definition \ref{def:critical}.  In fact $c_1 = c$ 
in this case, since both cells involve the same edges, the vertices in 
both are blocked, and each component of $T - E(c)$ contains the same 
number of vertices from each of $c_1$ and $c$.

Now suppose that $c_1$ has respectful edges.  Let $e$ be the
smallest such (recall that ``smallest" means that $\iota(e)$ is
minimal).  Since the edge $e$ is disrespectful in $c$, there is
some vertex $v \in c$ adjacent to $\tau(e)$ and satisfying
  $$ 0 < g( \tau(e), v) < g( \tau(e), \iota(e)).$$
Let $C$ be the connected component of $T - E(c)$
containing $\tau(e)$ and lying in the direction $g( \tau(e), v)$ from
$\tau(e)$.  This component contains vertices of $c$ and thus vertices of
$c_1$, since $c_1 \sim c$.  If $C$ contains unblocked vertices of $c_1$,
then any such vertex $v_1$ satisfies $v_1 < \iota(e)$, and so it follows
that $c_1$ is redundant.  If $c$ contains only blocked vertices, then it
follows that the vertex $v$ is a vertex of $c_1$, whence the edge $e$ is
disrespectful in $c_1$, a contradiction.
\end{proof}

\subsection{The complex $\widehat{\ud{n}{T}}$}

Define a complex $\widehat{\ud{n}{T}}$ as follows.  For each equivalence
class $[c]$ of $1$-cells in the set $K$ of open cells of $\ud{n}{T}$,
introduce a copy of $S^{1}$, denoted $S^{1}_{[c]}$.  Give $S^1_{[c]}$ a cell 
structure with one open $1$-cell, denoted $e^{1}_{[c]}$, and one $0$-cell.  
Form the finite product $\prod_{[c]} S^{1}_{[c]}$.  Since we are
interested in giving each cell an explicit characteristic map, we order
the factors in the product as follows.  Assign to each equivalence class
$[c]$ of $1$-cells the number $N([c])$ of the vertex $\iota(e)$, where $e$
is the unique edge satisfying $e \in c$.  This numbering of equivalence
classes is well-defined (though not one-to-one).  Arrange the factors of
$\prod_{[c]} S^{1}_{[c]}$ so that if $N([c_1]) < N([c_2])$, then the
factor $S^{1}_{[c_1]}$ occurs before $S^{1}_{[c_2]}$.  (This arrangement
of factors is not unique.)

Since each $1$-cell of this product corresponds naturally to an
equivalence class $[c]$ of $1$-cells in $\ud{n}{T}$, each $i$-cell
corresponds to a collection $\{ [c_1], \ldots, [c_i] \}$ of (distinct)
equivalence classes of $1$-cells.  We obtain the space
$\widehat{\ud{n}{T}}$ by throwing out an open $i$-cell $\{ [c_1], \ldots,
[c_i] \}$ if $\{ [c_1], \ldots, [c_i] \}$ has no upper bound.  If
$\{ [c_1], \ldots, [c_i] \}$ has an upper bound, then it has a least
upper bound $[c]$, and we label the corresponding $i$-cell by $[c]$.  
Note that, by Lemma \ref{lem:big}(3), the equivalence classes $[c]$ are in
one-to-one correspondence with cells of $\widehat{\ud{n}{T}}$, and the
dimension of the cell $c$ of $\ud{n}{T}$ is the same as that of the cell
labelled $[c]$ in $\widehat{\ud{n}{T}}$.

If $R$ is a field, then the \emph{exterior algebra} (\cite{Hatcher}, pg. 
217) on a set $\{ v_{1} , v_{2}, \ldots, v_{n+1} \}$, denoted 
$\Lambda_{R} [ v_{1}, \ldots, v_{n+1} ]$, is the $R$-module having 
relations $v_i v_j = - v_j v_i$ and $v_{i}^{2} = 0$.  The products 
$v_{i_1}v_{i_2} \ldots v_{i_j}$ ( $ 0 \leq j \leq n$; $i_1 < i_2 < 
\ldots < i_j$) form a basis.  The empty product is the multiplicative 
identity.

For any equivalence class $[c]$ of $j$-cells, let 
$\hat{\phi}_{[c]}:C_*(\hat{\ud{n}{T}}) \to \Z/2\Z$
denote the $j$-cocycle satisfying $\hat{\phi}_{[c]}([c]) = 1$ and
$\hat{\phi}_{[c]}([c']) = 0$ for all $[c'] \neq [c]$.  (Here $C_{j}(\hat{\ud{n}{T}})$
denotes the free abelian group generated by oriented $j$-cells of $\hat{\ud{n}{T}}$.  Note
that $\hat{\phi}_{[c]}$ is necessarily a cocycle, since all of the boundary maps
in the cellular chain complex of a subcomplex of $\prod S^{1}$ are $0$.)

\begin{proposition} \label{prop:udn}
The space $\widehat{\ud{n}{T}}$ is a CW complex.  The mod $2$ cohomology 
ring $H^{\ast} \left( \widehat{\ud{n}{T}} ; \mathbb{Z}/2\mathbb{Z} 
\right)$ is isomorphic to
  $$ \Lambda_{\mathbb{Z}/2\mathbb{Z}} \big[ [c_1], \ldots, [c_k] \big] 
  / I,$$
where $\{ [c_1], \ldots, [c_k] \}$ is the collection of all equivalence
classes of $1$-cells in $\ud{n}{T}$, and $I$ is the ideal generated by the
set of all monomial terms $[c_{i_1}]\cdot [c_{i_2}] \cdot \ldots \cdot [c_{i_m}]$ 
where $\{ [c_{i_1}],
\ldots, [c_{i_m}] \}$ has no upper bound.

The isomorphism sends the $j$-cocycle $\hat{\phi}_{[c]}$ to the product 
$[c_{i_1}]\cdot[c_{i_2}]\cdot\ldots\cdot[c_{i_j}]$,
where $\{ [c_{i_1}], \ldots, [c_{i_j}] \}$ is the unique collection with 
least upper bound $[c]$ as in Lemma \ref{lem:big}(3).

The elements $\hat{\phi}_{[c]}$ form a basis for the cohomology as $[c]$
ranges over all possible equivalence classes.

\end{proposition}

\begin{proof}

The complex $\widehat{\ud{n}{T}}$ inherits a cell structure from 
$\prod_{[c]} S^{1}_{[c]}$.  To prove that $\widehat{\ud{n}{T}}$ is a CW 
complex, we need to verify that, for every open cell in 
$\widehat{\ud{n}{T}}$, the attaching map
to $\prod_{[c]} S^{1}_{[c]}$ also maps into
$\widehat{\ud{n}{T}}$.  
This means we 
must show that if an open cell $\{ [c_{j_1}], \ldots, [c_{j_l}] \}$ is 
thrown out of the product $\prod_{[c]} S^{1}_{[c]}$, then so is any 
other open cell having $\{ [c_{j_1}], \ldots, [c_{j_l}] \}$ as a face.  
By definition, if the cell $\{ [c_{j_1}], \ldots, [c_{j_l}] \}$ is 
thrown out, the collection $\{ [c_{j_1}], \ldots, [c_{j_l}] \}$ has no 
upper bound.  Any cell having $\{ [c_{j_1}], \ldots, [c_{j_l}] \}$ as a 
face must be labelled by a collection $S$ of equivalence classes of 
$1$-cells satisfying $\{ [c_{j_1}], \ldots, [c_{j_l}] \} \subseteq S$.  
But as $\{ [c_{j_1}], \ldots, [c_{j_l}] \}$ has no upper bound, $S$ can 
have no upper bound.  Since $S$ has no upper bound, by the defininition 
of $\widehat{\ud{n}{T}}$ the open cell labelled $S$ is also thrown out 
of the product, as required.

The remaining statements then follow easily from the description of
$\widehat{\ud{n}{T}}$ as a subcomplex of $\prod_{[c]} S^{1}_{[c]}$, and
from the description in \cite{Hatcher} (pg. 227) of the cohomology rings
of subcomplexes of the torus. 
\end{proof}

\subsection{A map $q: \ud{n}{T} \rightarrow \widehat{\ud{n}{T}}$ and the
induced map on cohomology} 

For each edge $e$ in the tree $T$, choose a characteristic map $h_e :
[0,1] \rightarrow e$, such that $h_e (0) = \iota(e)$ and $h_e (1) =
\tau(e)$.  These maps induce characteristic maps on the cells of
$\ud{n}{T}$ as follows.  Let $c = \{ e_1, \ldots, e_i, v_{i+1}, \ldots,
v_n \}$ be an $i$-cell of $\ud{n}{T}$, and suppose without loss of
generality that $\iota(e_1) < \iota(e_2) < \ldots < \iota(e_i)$.  Define 
the
characteristic map $h_c : [0,1]^{i} \rightarrow c$ by 
  $$h_c (t_1, \ldots, t_i) = \{ h_{e_1}(t_1), \ldots, h_{e_i}(t_i), 
  v_{i+1},\ldots, v_{n} \}.$$
Note that $h_c (t_1, \ldots, t_i)$ is an $n$-element subset of $T$ here,
rather than an $n$-tuple of cells in $T$.  The map $h_c$ is a
homeomorphism.

We now choose characteristic maps for the cells of $\widehat{\ud{n}{T}}$.  
Begin by choosing a characteristic map $h_{[c]}: [0,1] \rightarrow
S^{1}_{[c]}$ for each $1$-cell $e^{1}_{[c]}$.  Suppose that $[c]$ is the
label of an $i$-dimensional cell in $\widehat{\ud{n}{T}}$.  Thus, $[c]$ is
the least upper bound of a collection $\{ [c_1], \ldots, [c_i] \}$ of
distinct equivalence classes of $1$-cells.  If $e_1, \ldots, e_i$ are the
unique edges satisfying $e_j \in c_j$ $(1 \leq j \leq i)$, then $e_1,
\ldots, e_i$ are pairwise disjoint by Lemma \ref{lem:big} (2).  In
particular, the natural numbers $N([c_1]), \ldots, N([c_i])$ are all
different.  We assume, without loss of generality, that $N([c_1]) < \ldots
< N([c_i])$.  The characteristic map for the cell $[c]$ in
$\widehat{\ud{n}{T}}$ is $\hat{h}_{[c]}: [0,1]^{i} \rightarrow
e^{i}_{[c]}$, defined by $\hat{h}_{[c]} (t_1, \ldots, t_i) = (
h_{[c_1]}(t_1), \ldots, h_{[c_i]}(t_i) )$, where the value of
$\hat{h}_{[c]}(t_1, \ldots, t_i)$ on each omitted factor $[\hat{c}]$, where
$[\hat{c}]$ is an equivalence class of $1$-cells, is always assumed to be
the unique vertex of $S^{1}_{[\hat{c}]}$.

Now we are ready to define the map $q: \ud{n}{T} \rightarrow
\widehat{\ud{n}{T}}$.  Consider the following diagram:

\begin{figure}[!h]
\[
\xymatrix{
  \coprod_{c \in K} [0,1]^{(\dim c)} 
    \ar[d]_{\coprod_{c\in K} h_c}
    \ar[dr]^{\coprod_{c\in K} \hat{h}_{[c]}} & \\
  \ud{n}{T} \ar@{-->}[r]_{q} &
  \widehat{\ud{n}{T}}
}
\]
\caption{Defining the map $q$.}
\label{fig:diagram}
\end{figure}

The vertical arrow is a quotient map, so, by a well-known principle (e.g.,
\cite{Dugundji}, Theorem 3.2 ), there will exist a well-defined map $q$ making
the above diagram commute if $\coprod_{c \in K} \hat{h}_{[c]}$ is constant
on point inverses of $\coprod_{c \in K} h_c$.

\begin{proposition} \label{prop:inducedmap}
There is a well-defined map $q: \ud{n}{T} \rightarrow \widehat{\ud{n}{T}}$
making the above diagram commute. The map $q^{\ast}: H^{\ast} \left(
\widehat{\ud{n}{T}}; \mathbb{Z}/2\mathbb{Z} \right) \rightarrow H^{\ast}
\left( \ud{n}{T} ; \mathbb{Z}/2\mathbb{Z} \right)$ sends
$\hat{\phi}_{[c]}$ to the cohomology class of $\phi_{[c]} \in C^{\ast}
\left( \ud{n}{T}; \mathbb{Z}/2\mathbb{Z} \right)$, where $\phi_{[c]}$ is
the cellular cocycle satisfying:
\begin{eqnarray*}
\phi_{[c]} \left( \tilde{c} \right) = 1 & & ~if~ \tilde{c} \sim c \\
\phi_{[c]} \left( \tilde{c} \right) = 0 & & ~otherwise~.
\end{eqnarray*}
\end{proposition}

\begin{proof}
We have to show that $\coprod_{c \in K} \hat{h}_{[c]}$ is constant on point
inverses of $\coprod_{c \in K} h_c$.  For this, it is sufficient to show
that $\hat{h}_{[c']} \circ h^{-1}_{c'} = \hat{h}_{[c]} \circ h^{-1}_{c}
\mid_{c'}$, where $c'$ is a codimension-one face of $c$.  We let $c$ be a
$j$-dimensional cell in $\ud{n}{T}$, say $c = \{ e_1, \ldots, e_j,
v_{j+1}, \ldots, v_{n} \}$, where the edges $e_i$ are arranged in order.  
Assume, without loss of generality, that $c' = \{ \iota(e_1), e_2, \ldots,
e_j, v_{j+1}, \ldots, v_n \}$.

Choose a point $x \in c'$. (Here we really mean a point in the cell $c'$,
rather than one of the ``members" $\iota(e_1), e_2, \ldots, e_j, v_{j+1},
\ldots, v_n$ of $c'$.)  Suppose that $x = h_{c'}(t_1, \ldots, t_{j-1})$,
i.e., $x = \{ \iota(e_1), h_{e_2}(t_1), \ldots, h_{e_j}(t_{j-1}), v_{j+1},
\ldots, v_{n} \}$.

\begin{eqnarray*}
\hat{h}_{[c']} \circ h^{-1}_{c'} (x) & = & \hat{h}_{[c']}(t_1, \ldots,t_{j-1}) \\ 
& = & \left( h_{[c'_{1}]}(t_1), \ldots, h_{[c_{j-1}]}(t_{j-1}) \right).\\
\end{eqnarray*}
\begin{eqnarray*}
\hat{h}_{[c]} \circ h^{-1}_{c}(x) & = & \hat{h}_{[c]}(0, t_1, \ldots, t_{j-1}) \\
& = & \left( h_{[c_1]}(0), h_{[c_2]}(t_1), \ldots, h_{[c_j]}(t_{j-1}) \right) \\
& = & \left( h_{[c_2]}(t_1), \ldots, h_{[c_j]}(t_{j-1}) \right).  
\end{eqnarray*}
For this last equality, recall that $h_{[c_1]}(0)$ is the vertex of
$S^{1}_{[c_1]}$, and we omit such factors for the sake of simplicity.

Now $[c'_{i}]$, by definition, is the unique equivalence class of
$1$-cells satisfying: (i) $e_{i+1} \in c'_{i}$, and (ii) $[c'_{i}] \leq
[c]$.  Note that $[c_{i+1}]$ has the same properties, so $[c'_{i}] =
[c_{i+1}]$.  It follows that the map $q$ exists.

The remaining statements about cohomology follow easily from the
description of the map $q: \ud{n}{T} \rightarrow \widehat{\ud{n}{T}}$.  
The main point is that the interior of a cell $c$ in $\ud{n}{T}$ is mapped
homeomorphically to the interior of $[c]$, and thus the mod $2$ mapping
number of $c$ with $[c]$ is equal to $1$.
\end{proof}

To describe the mod $2$ cohomology ring of $\ud{n}{T}$, we will need to
recall a result from \cite{Farley2}:

\begin{theorem} \emph{(}\cite{Farley2}, Theorem 3.7 \emph{)}
The boundary maps in the Morse complex $( M_{\ast}(\ud{n}{T}),
\tilde{\partial} )$ are all zero.  In particular, $H_{i}( \ud{n}{T} )$ is
a free abelian group of rank equal to the number of critical $i$-cells in
$\ud{n}{T}$.
\qed
\end{theorem}

By Lemma \ref{lem:homology}, we obtain explicit cycles in $C_{i}(
\ud{n}{T})$, the $i$th cellular chain group of $\ud{n}{T}$, by applying
the map $f^{\infty}$ to any linear combination of critical $i$-cells.  A
collection of representatives for a distinguished basis of the cellular
$i$-dimensional homology is thus $\{ f^{\infty}(c) \mid \, \, $c$
~is~a~critical~ i-cell \}$.  For the sake of simplicity in notation, we
express a cellular homology class in $H_{\ast}(\ud{n}{T})$ as a linear
combination $\Sigma a_i c_i$ of critical cells $c_i$, as opposed to
$f^{\infty} ( \Sigma a_i c_i )$.  Identify $H^{i}( \ud{n}{T};
\mathbb{Z}/2\mathbb{Z} )$ with $Hom( H_{i}( \ud{n}{T});
\mathbb{Z}/2\mathbb{Z} )$ by the universal coefficient isomorphism.  Let
$c^{\ast}$ denote the dual of a critical cell $c$.

\begin{proposition} \label{prop:cup} ~\newline   
\begin{enumerate}
\item If $c$ is a critical cell in $\ud{n}{T}$, then $q^{\ast}(
\hat{\phi}_{[c]} ) = c^{\ast}$.

\item Let $c$ be a critical cell in $\ud{n}{T}$.  If $[c]$ is the least
upper bound of $\{ [c_1], \ldots, [c_i] \}$, where the $[c_1], \ldots,
[c_i]$ are distinct equivalence classes of $1$-cells, 
then, without loss of generality, $c_1, \ldots, c_i$
are critical and
  $$ c^{\ast}_{1} \cup \ldots \cup c^{\ast}_{i} = c^{\ast}.$$

\item If $[c_1], \ldots, [c_i]$ are distinct equivalence classes of
$1$-cells having the least upper bound $[c]$, then
  $$ \left[ \phi_{[c_1]} \right] \cup \ldots \cup \left[ 
  \phi_{[c_i]} \right] = \left[ \phi_{[c]} \right].$$
If $[c_1], \ldots, [c_i]$ have no upper bound, then the above cup product
is $0$.
\end{enumerate}
\end{proposition}

\begin{proof}
(1)  By Proposition \ref{prop:inducedmap}, $\phi_{[c]}$ is a cocycle
representative of $q^{\ast}(\hat{\phi}_{[c]})$, where
$\phi_{[c]}(\tilde{c}) = 1$ if $\tilde{c} \sim c$, and
$\phi_{[c]}(\tilde{c}) = 0$ otherwise. Note that, by Lemma
\ref{lem:big}(5), the support of $\phi_{[c]}$ consists of redundant
cells, and a single critical cell ($c$ itself), but no collapsible cells.

We evaluate the cohomology class of $\phi_{[c]}$ on a basis for
$H_{i}(\ud{n}{T})$ consisting of critical cells $c_1, \ldots, c_j$.  As 
$c$ is the unique critical cell in the support of $\phi_{[c]}$, for a 
critical $1$-cell $c_k$ (viewed as a homology class), $\phi_{[c]}(c_k) = 
1$ if and only if $c_k = c$.  The statement of (1) follows.

(2) The statement that $c_1, \ldots, c_i$ may be chosen to be critical 
follows from Lemma \ref{lem:big}(4). By Proposition \ref{prop:udn}, 
$\hat{\phi}_{[c]} = \hat{\phi}_{[c_1]} \cup \ldots \cup 
\hat{\phi}_{[c_i]}$ in $H^{\ast}(\widehat{\ud{n}{T}}; 
\mathbb{Z}/2\mathbb{Z} )$.  By (1), the statement of (2) follows.

(3)  This follows from applying Propositions \ref{prop:udn} and
\ref{prop:inducedmap}.
\end{proof}

\subsection{A computation of $H^{\ast}( \ud{4}{T_{min}};
\mathbb{Z}/2\mathbb{Z})$}
\label{sec:computation}

\begin{figure}[!h]
\centering
\input{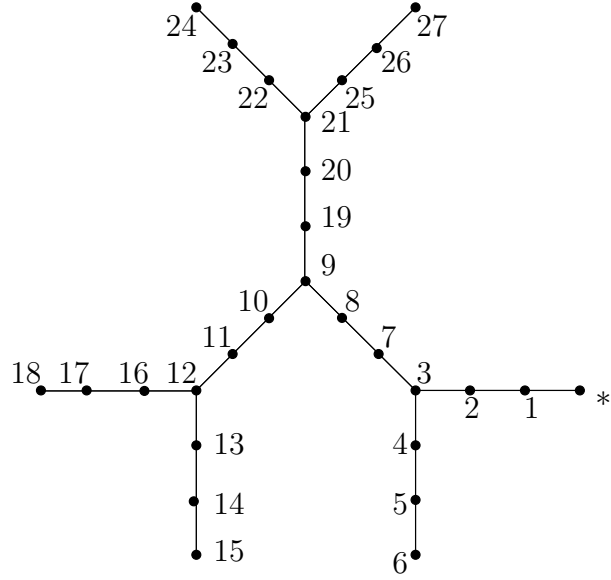}

\caption{The minimal nonlinear tree $T_{min}$.}
\label{fig:favtree}  
\end{figure}

Let $T_{min}$ be the tree depicted in Figure \ref{fig:favtree}.  The 
tree $T_{min}$ is the tree, unique up to homeomorphism, with the fewest 
number of essential vertices which is not `linear': i.e., the vertices 
of degree 3 or more in $T_{min}$ do not all lie on a single embedded 
line segment.  We compute the mod $2$ cohomology ring of 
$\ud{4}{T_{min}}$ as an application of the ideas of this section.  The 
results will be used in the proof of Theorem \ref{thm:Ghrist}.

To begin, we will need to compute the integral homology groups of
$\ud{4}{T_{min}}$.  According to Proposition 4.1 of \cite{Farley2}, we
have
\begin{eqnarray*}
H_0 (\ud{4}{T_{min}}) \cong \mathbb{Z} & , & H_1 (\ud{4}{T_{min}}) \cong\mathbb{Z}^{24} \\
H_2 (\ud{4}{T_{min}}) \cong \mathbb{Z}^{6} & , & H_n (\ud{4}{T_{min}}) \cong 0 (n \geq 3).    
\end{eqnarray*}

We also need to describe the critical cells $c$, which are determined by
the following choices.  First, choose the locations of the edges of $c$.  
The requirement that the edges in $c$ be disrespectful means that
there are only four possibilities:  $[v_3 , v_7], [v_9 , v_{19}], [v_{12},
v_{16}]$, and $[v_{21}, v_{25}]$.  We let $e_k$ denote the unique edge in
$T$ having $v_{k}$ as its initial vertex.  With this notation, we can
rewrite the four edges above as $e_7$, $e_{19}$, $e_{16}$, and $e_{25}$,
respectively.  A choice of $n$ edges from this collection, together with
the requirement that the edges be disrespectful, forces $n$
vertices also to be in $c$.  More specifically, $v_4 \in c$ if $e_7 \in
c$, $v_{10} \in c$ if $ e_{19} \in c$, $v_{13} \in c$ if $e_{16} \in c$,
and $v_{22} \in c$ if $ e_{25} \in c$.  It immediately follows that there
are exactly $\binom{4}{2} = 6$ critical $2$-cells and no critical
$n$-cells for $n \geq 3$, so $H_{2}( \ud{4}{T_{min}} ) \cong
\mathbb{Z}^{6}$ and $H_{3}( \ud{4}{T_{min}} ) \cong 0$ if $n \geq 3$.

If $c$ is a critical $1$-cell, there is another choice to make.  We have
determined the location of an edge $e$ and a vertex $v$.  The locations of
the other two vertices are completely determined by specifying how many
are in each component of $T - \{ \tau(e) \}$, since all vertices in $c$
are blocked.  There are $3$ distinguishable components of $T - \{ \tau(e)
\}$ and $2$ remaining indistinguishable vertices, which make $6$ possible
ways.  There are thus a total of $24$ critical $1$-cells $c$, since there
are $4$ possible choices for the edge $e$ in $c$, $6$ choices for the
remaining vertices and these choices are independent.

Finally, we note that there is only one critical $0$-cell.

To understand the multiplication in $H^{\ast}( \ud{4}{T_{min}};
\mathbb{Z}/2\mathbb{Z} )$, it is clearly enough to understand the product
of any two elements $c_{1}^{\ast}$, $c_{2}^{\ast}$ ($c_1$ and $c_2$ are
critical $1$-cells) of the standard dual basis for $H^{1}(
\ud{4}{T_{min}}; \mathbb{Z}/2\mathbb{Z})$.  Proposition \ref{prop:cup}(1)
says that $c_{1}^{\ast} = \left[\phi_{[c_1]} \right]$ and
$c_{2}^{\ast} = \left[\phi_{[c_2]} \right]$, and Proposition
\ref{prop:cup}(3) says that $\left[ \phi_{[c_1]} \right] \cup \left[
\phi_{[c_2]} \right] =0$ if $c_1 \sim c_2$, or if $\{ [c_1], [c_2]
\}$ has no upper bound.  We are therefore led to determine which
equivalence classes $[c]$ of $2$-cells can be the upper bound for a pair
of distinct equivalence classes of $1$-cells $[c_1]$ and $[c_2]$, where
$c_1$ and $c_2$ are critical.

For this, it will be helpful to have a definition.  If $c$ is a $j$-cell
in $\ud{n}{\Gamma}$, and $c'$ is obtained from $c$ by replacing each
member of some collection $\{ e_{i_1}, \ldots, e_{i_m} \} \subseteq E(c)$
of edges with either its initial or its terminal vertex, then $c'$ is the
result of \emph{breaking} the edges $\{ e_{i_1}, \ldots, e_{i_m} \}$ in
$c$. We note that the choice of replacing a given edge $e_{i_l}$ with
$\tau ( e_{i_l} )$ or $\iota ( e_{i_l} )$ is made independently for each
edge, and these choices do not affect $[c']$.  In fact, by Lemma
\ref{lem:big} (1), $[c']$ is completely determined by its edges, since
$[c'] \leq [c]$.  Note also that $c'$ is a codimension-$m$ face of $c$,
and conversely if $c'$ is a codimension-$m$ face of $c$ then $c'$ is
obtained from $c$ by breaking $m$ edges.

Suppose that $[c]$ is the least upper bound for two distinct equivalence 
classes $[c_1]$ and $[c_2]$, where $c_1$ and $c_2$ are critical 
$1$-cells in $\ud{4}{T_{min}}$.  For $i=1,2$, let $e_i$ be the unique 
edge in $c_i$. (Note: this is inconsistent with the convention that 
$\iota(e_k) = v_k$, but should cause no confusion.)  By Lemma 
\ref{lem:big}(2), $e_1$ and $e_2$ are disjoint edges, and, by the above 
description of critical cells in $\ud{4}{T_{min}}$, $\{ e_1, e_2 \} 
\subseteq \{ e_7 , e_{19}, e_{16}, e_{25} \}$.  Thus, if $[c]$ is an 
upper bound for $[c_1]$ and $[c_2]$ where $c_1$ and $c_2$ critical, we 
have: (i) the edges of $c$, namely $e_1$ and $e_2$, are distinct 
elements of the above $4$-element set, and (ii) the cell resulting from 
breaking either of the edges $e_1, e_2$ in $c$ must be equivalent to a 
critical cell.  Reformulating slightly, we get the following conditions, 
which must be satisfied by $c$:
\begin{enumerate}
\item $\{ e_1, e_2 \} \subseteq \{ e_7 , e_{19} , e_{16} , e_{25} \}$.

\item If $e_7 \in c$, then $v_4$, $v_5$, or $v_6$ is in $c$.

\item If $e_{16} \in c$, then $v_{13}$, $v_{14}$, or $v_{15}$
is in $c$.

\item If $e_{25} \in c$, then $v_{22}$, $v_{23}$, or $v_{24}$
is in $c$.

\item If $e_{19}  \in c$, then either $e_{16} \in c$
or at least one element of $\{ v_{10}, \ldots, v_{18} \}$ is in $c$.
\end{enumerate} 
In light of (3), (5) may be replaced with:
\begin{enumerate}[]
\item{(5')} If $e_{19} \in c$, then at least one element of $\{ 
v_{10}, \ldots, v_{18} \}$ is in $c$.
\end{enumerate}

A total of $10$ distinct equivalence classes of $2$-cells have
representatives $c$ satisfying (1)-(5).  Representatives of these classes
are:

\begin{figure}[!t]
\centering
\input{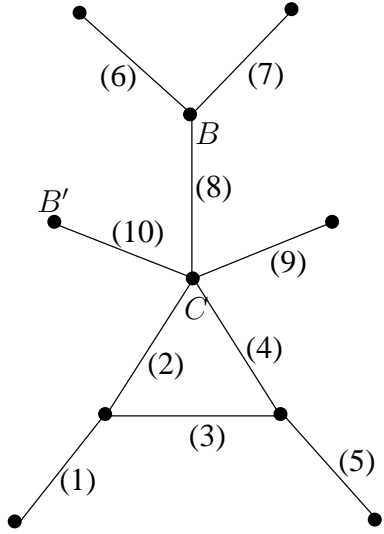}

\caption{A preliminary picture of the relations in $H^{\ast}\left( 
\ud{n}{T_{min}} ; \mathbb{Z} / 2\mathbb{Z} \right)$. The numbers on the 
edges refer to the numbering of $2$-cells in the text.  The vertices 
represent duals of critical $1$-cells.  Here, $B$ denotes the critical 
$1$-cell $\{e_{19}, v_{12}, v_{11}, v_{10}\}$, $B'$ denotes the critical 
$1$-cell $\{e_{19}, v_{11}, v_{10}, \ast\}$, and $C$ denotes the 
critical $1$-cell $\{e_{16}, v_{13}, v_1, \ast\}$.}
\label{fig:prod}
\end{figure}

\begin{center}
\begin{tabular}{l l}
$(1) \, \, \{ e_7 , e_{19}, v_{4}, v_{10} \}$ & $\quad (6) \, \, \, \{ 
e_{19}, e_{16}, v_{13}, v_{14} \}$ \\
$(2) \, \, \{ e_7 , e_{16}, v_{4}, v_{13} \}$ & $\quad (7) \, \, \, \{ 
e_{19}, e_{16}, v_{13}, v_{17} \}$ \\
$(3) \, \, \{ e_7 , e_{25}, v_{4}, v_{22} \}$ & $\quad (8) \, \, \,  \{ 
e_{19}, e_{16}, v_{13}, v_{10} \}$ \\
$(4) \, \, \{ e_{16} , e_{25}, v_{13}, v_{22} \}$ & $\quad (9) \, \, \, \{ 
e_{19}, e_{16}, v_{13}, v_{20} \}$ \\
$(5) \, \, \{ e_{19} , e_{25}, v_{10}, v_{22} \}$ & $\quad (10) \, \, \{ 
e_{19}, e_{16}, v_{13}, \ast \}$  
\end{tabular}
\end{center}
 
\vspace{10pt}

The reason is that a choice of edges $e_1$, $e_2$ from $\{ e_7 , e_{16},
e_{19}, e_{25} \}$ completely determines $[c]$, by (1)-(5), unless $\{
e_1, e_2 \} = \{ e_{16}, e_{19} \}$.  The first five edges listed above
result from the five cases in which $\{ e_1, e_2 \} \neq \{ e_{16}, e_{19}
\}$.  If $\{ e_1, e_2 \} = \{ e_{16}, e_{19} \}$, then one of the vertices
of $c$ must be $v_{13}$, $v_{14}$ or $v_{15}$ (by (3)), but the other
vertex may be chosen from any of the five remaining components of $T_{min}
- ( e_{16} \cup e_{19} )$, and this accounts for the last five $2$-cells
above.

Figure \ref{fig:prod} below depicts these $10$ equivalence classes $[c]$
of $2$-cells as line segments whose endpoints are the equivalence classes
$[c_1], [c_2]$ of $1$-cells ($c_1, c_2$ critical) satisfying $[c_1], [c_2]
\leq [c]$.

It follows from what we've said so far that there are only $10$
equivalence classes $[c']$ of $1$-cells such that: (1) $[c']$ contains a
critical $1$-cell, and (2) there is another distinct equivalence class
$[c'']$, also containing a critical $1$-cell, such that $\{ [c'], [c'']
\}$ has an upper bound.  By Proposition \ref{prop:cup}(1) \& (3), these
equivalence classes correspond to the only elements of the standard basis
for $H^{1}( \ud{4}{T_{min}}; \mathbb{Z}/2\mathbb{Z})$ which might have
non-trivial cup products.  In fact, checking labels of edges, it is not
difficult to see that cases (1)-(5) and (8) are all (distinct)  critical
cells, and thus correspond to linearly independent elements of the
standard basis for $H^{2}(\ud{4}{T_{min}}; \mathbb{Z}/2\mathbb{Z})$.  We
now consider the remaining edges.

\begin{figure}[!b]
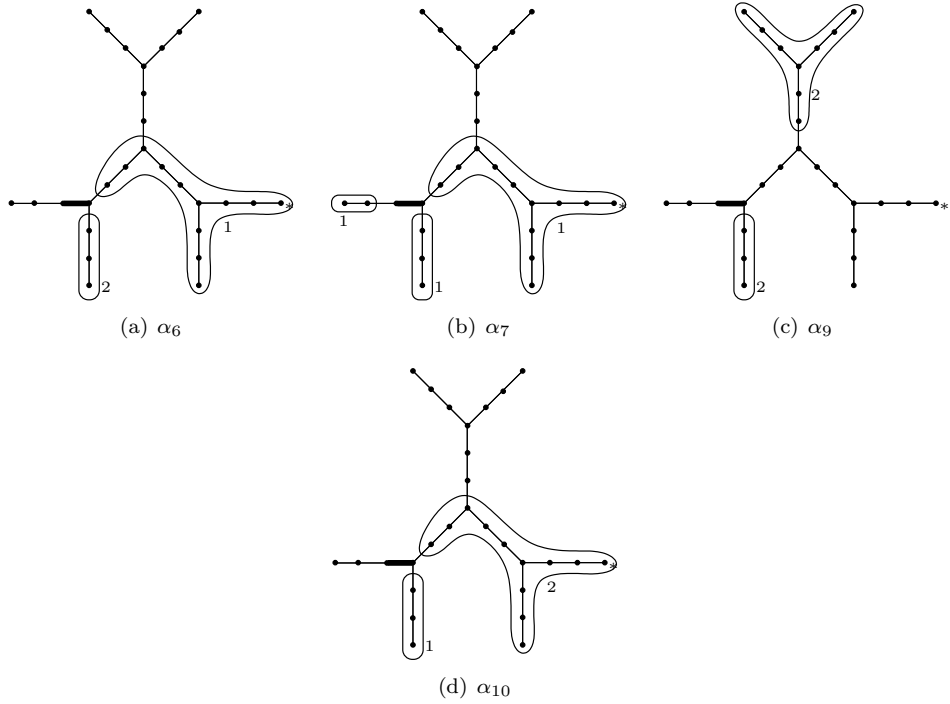

\subfigure[$\alpha_6$]{
  \input{a6.pstex_t}
}\hfill
\subfigure[$\alpha_7$]{
  \input{a7.pstex_t}
}\hfill
\subfigure[$\alpha_9$]{
  \input{a9.pstex_t}
}\hfill
\subfigure[$\alpha_{10}$]{
  \input{a10.pstex_t}
}

\caption{The $1$-cochains $\alpha_6$, $\alpha_7$, $\alpha_9$, and 
$\alpha_10$ used in Lemmas \ref{lem:alpha1} and \ref{lem:alpha2}.  Each 
$1$-cochain maps a given $1$-cell $c$ to $1$ in $\mathbb{Z}/2\mathbb{Z}$ 
if and only if $c$ contains the edge shown and has exactly as many 
strands specified in each circled portion of the tree.}
\label{fig:alphas}
\end{figure}

\begin{lemma}\label{lem:alpha1}
If $\hat{c}$ is one of the $2$-cells (6), (7) or (9), then
$\phi_{[\hat{c}]}$ represents $0$ in cohomology.
\end{lemma}

\begin{proof}
Consider the following cochains $\alpha_6, \alpha_7, \alpha_9: 
C_{1}(\ud{4}{T_{min}}) \rightarrow \mathbb{Z}/2\mathbb{Z}$, depicted 
pictorially in Figure \ref{fig:alphas}. The cochain 
$\alpha_6$ sends a given $1$-cell $c$ to $1$ if and only if: (1) $c$ 
contains the edge $e_{16}$; (2) $c$ contains exactly one element from 
$\{ \ast, v_{1}, \ldots, v_{11} \}$, and (3) $c$ contains exactly two 
elements from $\{ v_{13}, v_{14}, v_{15} \}$.  (Of course, $\alpha_6$ 
sends any other $1$-cell to $0$.)  The cochain $\alpha_7$ sends a given 
$1$-cell $c$ to $1$ if and only if $c$ satisfies conditions (1) and (2) 
in the definition of $\alpha_6$, as well as: (3') $c$ contains exactly 
one element from $\{ v_{17}, v_{18} \}$, and exactly one element from 
$\{ v_{13}, v_{14}, v_{15} \}$.  The cochain $\alpha_9$ sends a given 
$1$-cell $c$ to $1$ if and only if $c$ contains exactly two elements 
from the set $\{ v_{19}, \ldots, v_{27} \}$, exactly one element from 
the set $\{ v_{13}, v_{14}, v_{15} \}$, and the edge $e_{16}$.

We leave it as an exercise to show that the coboundaries 
$\delta(\alpha_6)$, $\delta(\alpha_7)$, and $\delta(\alpha_9)$, are 
precisely $\phi_{[\hat{c}]}$, where $\hat{c}$ is the cell (6), (7), and 
(9), respectively.
\end{proof}

\begin{lemma}\label{lem:alpha2}
Let $c_1$ be the cell labelled (8), and let $c_2$ be the cell labelled
(10).  The cocycles $\phi_{[c_1]}$ and $\phi_{[c_2]}$ are cohomologous in 
$H^{1}( \ud{4}{T_{min}}; \mathbb{Z}/2\mathbb{Z})$.
\end{lemma}

\begin{proof}
Consider the $1$-cochain $\alpha_{10}$ defined as follows and depicted 
in Figure \ref{fig:alphas}.  Let $\alpha_{10}$ be the $1$-cochain which 
sends a $1$-cell $c$ to $1$ if and only if: (1) $c$ contains the edge 
$e_{16}$; (2) $c$ contains exactly one of the vertices $\{ v_{13}, 
v_{14}, v_{15} \}$, and (3) $c$ contains exactly two vertices from $\{ 
\ast, v_{1}, \ldots, v_{11} \}$.

We leave it as an exercise to show that $\delta(\alpha_{10}) = 
\phi_{[c_1]} + \phi_{[c_2]}$.
\end{proof}

\begin{figure}[!h]
\centering
\input{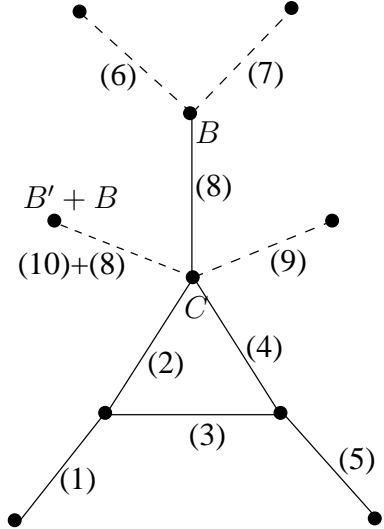}

\caption{A picture of the cohomology ring $H^{\ast} \left( \ud{n}{T_{min}}
; \mathbb{Z}/ 2\mathbb{Z} \right)$.  The solid edges represent duals of
critical $2$-cells;  the vertices represent (sums of) duals of critical 
$1$-cells. Two vertices cup to the solid edge connecting them, or to $0$ 
if there is no such edge.}
\label{fig:prodtwo}
\end{figure}

We now interpret Figure \ref{fig:prod} as a multiplication table for the
cup product.  If we let dashed edges correspond to $0$ products and
perform an elementary row operation, we arrive at Figure
\ref{fig:prodtwo}.

We thus arrive at a complete description of the multiplication in 
$H^{\ast}( \ud{4}{T_{min}}; \mathbb{Z}/2\mathbb{Z})$: an 
$18$-dimensional subspace $W$ of $H^{1} ( \ud{4}{T_{min}} ; \mathbb{Z}/ 
2\mathbb{Z} )$ annihilates all one-dimensional cohomology classes. This 
subspace $W$ is spanned by the duals of the $14$ critical $1$-cells not 
appearing in Figure \ref{fig:prodtwo}, together with the four elements 
of $H^{1} \left( \ud{n}{T} ; \mathbb{Z} / 2\mathbb{Z} \right)$ which are 
endpoints of only dashed lines.  The multiplication in the remaining 
$6$-dimensional subspace is described by the subgraph of Figure 
\ref{fig:prodtwo} consisting of six solid lines and the six vertices 
they connect: two basis elements cup to the label of the solid edge 
connecting them, or to $0$ if there is no such edge.

This description of the multiplication in $H^{\ast} \left(
\ud{4}{T_{min}}; \mathbb{Z}/2\mathbb{Z} \right)$ suggests that it is an
exterior face algebra, an idea we define in the next section.

\section{Counterexamples to Ghrist's Conjecture}\label{sec:Counterex}

\subsection{Preliminaries on Exterior Face Rings}

For a ring $R$, an \emph{exterior ring} over $R$ on a set $\{ v_{1} ,
v_{2}, \ldots, v_{n} \}$, denoted $\Lambda_{R} [ v_{1}, \ldots, v_{n}
]$, is the free $R$-module having the products $v_{i_1}v_{i_2} \ldots
v_{i_j}$ ($0 \leq i_1 < i_2 < \ldots < i_j \leq n$) as a basis.  
The empty product is the multiplicative identity.  The multiplicative
relations are generated by all relations of the following types:  $v_i v_j
= - v_j v_i$ and $v_{i}^{2} = 0$.

Let $K = \left( \{ v_1, \ldots, v_n \}, S \right)$ be a finite simplicial
complex and $R$ be a commutative ring with identity. The \emph{exterior
face ring} $\Lambda_{R} \left( K \right)$ of $K$ over $R$ is the quotient
of the exterior ring $\Lambda_{R}[v_1, \dots, v_n]$ by the relations
$v_{i_1} \ldots v_{i_k} = 0$ for $0 \leq i_1 < i_2 < \ldots < i_k \leq n$
and when $\{ v_{i_1}, \ldots, v_{i_k} \} \not \in S$. Note that an
exterior ring on a set $\{ v_{1} , v_{2}, \ldots, v_{n}, v_{n+1} \}$ is
isomorphic the exterior face ring $\Lambda_{R} \left( K \right)$ where $K$
is a standard $n$-simplex.  If $R$ is a field, then $\Lambda_R(K)$
inherits an algebra structure, and is called an \emph{exterior face
algebra}.

\begin{example}  
The calculation of the previous subsection shows that the ring
$H^{\ast}(\ud{4}{T_{min}}; \mathbb{Z}/2\mathbb{Z})$ is isomorphic to
$\Lambda_{\mathbb{Z}/2\mathbb{Z}}(K)$, where $K$ is the union of $18$
isolated vertices with a graph isomorphic to the one in Figure
\ref{fig:prodtwo} (2) consisting of the six solid edges labelled (1) - (5)
and (8), and vertices incident with them.
\end{example}

If $R = \mathbb{Z}/2\mathbb{Z}$, then the exterior face ring
$\Lambda_{\mathbb{Z}/2\mathbb{Z}}\left( K \right)$ is a quotient of a
polynomial ring:
  $$ \Lambda_{\mathbb{Z}/2\mathbb{Z}} \left( K \right) = 
  \mathbb{Z}/2\mathbb{Z} \left[ v_1, \ldots, v_n \right]/ I(K),$$
where $I(K)$ is the ideal of $\mathbb{Z}/2\mathbb{Z} \left[ v_1, \ldots, 
v_n \right]$ generated by the set
  $$ \left\{ v_{1}^{2}, \ldots, v_{n}^{2} \right\}
  \cup \left\{ v_{i_1} v_{i_2} \ldots v_{i_k} \mid 
  \left\{ v_{i_1}, \ldots, v_{i_k} \right\} \not \in S \right\}.$$
In this case, since $R$ is a field, we have that $\Lambda_R(K)$
inherits an algebra structure. Throughout the rest of the paper, all
exterior face rings $\Lambda_{R}(K)$ will be over
$\mathbb{Z}/2\mathbb{Z}$, and we will therefore drop the subscript $R$
without further comment.

A simplicial complex $K$ is \emph{flag} if, whenever a collection of
vertices $v_{i_1}, \ldots, v_{i_j} \in K$ pairwise span edges, $\{
v_{i_1}, \ldots, v_{i_j} \}$ is a simplex of $K$.

In case $K$ is a flag complex, there is a simple set of generators for
$I(K)$:

\begin{lemma} \label{flaglemma} 
If $K$ is a flag complex, then 
  $$ I(K) = 
  \left\langle \left\{ v_{1}^{2}, \ldots, v_{n}^{2} \right\}
  \cup \left\{ v_{i_1}v_{i_2} \mid \left\{ v_{i_1}, v_{i_2} \right\} 
  \not \in S \right\} \right\rangle.$$
\end{lemma}

\begin{proof}
Let $I'(K)$ denote the ideal on the right half of the equality in the
lemma. We need to show that $I(K) \subseteq I'(K)$, the reverse inclusion
being obvious.

Suppose that $v_{i_1} \ldots v_{i_k}$ satisfies $\left\{ v_{i_1}, \ldots,
v_{i_k} \right\} \not \in S$.  Since $K$ is a flag complex, there must
exist $v_{j_1}, v_{j_2} \in \left\{ v_{i_1}, \ldots, v_{i_k} \right\}$
such that $\left\{ v_{j_1}, v_{j_2} \right\} \not \in S$.  It follows that
$v_{j_1} v_{j_2} \in I'(K)$.  Now $v_{j_1}v_{j_2} \mid v_{i_1}\ldots
v_{i_k}$, so $v_{i_1} \ldots v_{i_k} \in I'(K)$.  Thus, $I(K) \subseteq
I'(K)$. 
\end{proof}

\begin{example} \label{ex:RAAG} Let $\Gamma$ be a finite simple graph.  
The \emph{right-angled Artin group} $G_{\Gamma}$ associated to $\Gamma$ 
is a group defined by a presentation in which the generators are in 
one-to-one correspondence with vertices of $\Gamma$, and relations 
consist of all commutators of the form $[v_i , v_j ]$, where $v_i$ and 
$v_j$ are adjacent in $\Gamma$.

Charney and Davis \cite{CharneyDavis} have described $K(G_{\Gamma},1)$ 
complexes for all right-angled Artin groups (generalizing the Salvetti 
complex of \cite{Salvetti} for spherical Artin groups).  Begin with a 
torus $\prod S^{1}$, where the factors are in one-to-one correspondence 
with vertices in $\Gamma$.  Assume that each $S^{1}$ is given the 
standard cellulation, consisting of one $0$-cell and one $1$-cell.  
Their $K(G_{\Gamma}, 1)$ space is obtained from this product by throwing 
out an open $i$-cell if the $i$ $1$-cells in its factorization 
correspond to vertices $v_1 , v_2 , \ldots , v_i$ which do not form a 
clique, i.e., if some pair of vertices $v_{j_1}, v_{j_2} \in \{ v_1, 
\ldots, v_i \}$ do not span an edge of $\Gamma$.

This description of $K(G_{\Gamma},1)$, together with the description of
the cohomology rings of subcomplexes of a torus in \cite{Hatcher} (pg.
227), implies 

\begin{proposition} \label{prop:flag} The cohomology ring $H^{\ast}(
G_{\Gamma}; \mathbb{Z} / 2\mathbb{Z} )$ is the exterior face ring $\Lambda
(K)$, where $K$ is the unique flag complex having $\Gamma$ as its
$1$-skeleton. \qed
\end{proposition}

\noindent (Here $K$ is the simplicial complex whose $n$-simplices are the 
cliques in $\Gamma$ having $n+1$ members.)
\end{example}

We can now give a simple principle which will allow us to find
counterexamples to Ghrist's Conjecture \ref{conj:Ghrist}.  A homomorphism
$\phi : R [x_1 , \ldots , x_l ] \rightarrow R [y_1 , \ldots , y_m ]$
between polynomial rings is \emph{degree-preserving} if it sends any
homogeneous polynomial of degree $k$ to another homogeneous polynomial of
degree $k$ (or, equivalently, if it sends any homogeneous polynomial of
degree $1$ to another homogeneous polynomial of degree $1$).  More
generally, if $R_1$ and $R_2$ are quotients of polynomial rings by ideals
generated by homogeneous polynomials, then $\phi : R_1 \rightarrow R_2$ is
\emph{degree-preserving} if any equivalence class of homogeneous
polynomials is mapped to an equivalence class of homogeneous polynomials
of the same degree.

\begin{proposition} \label{prop:principle} 
Let $K$ be a flag complex, and let $\partial \Delta^{n}$ be the boundary
of the standard $n$-simplex ($n \geq 2$). If $\phi : \Lambda (K)
\rightarrow \Lambda \left( \partial \Delta^{n} \right)$ is a
degree-preserving surjection, then $\ker\phi$ cannot be generated by
homogeneous degree $1$ and degree $2$ elements.
\end{proposition}

\begin{proof}
Let $\left\{ v_1, \ldots, v_m \right\}$ be the vertices of $K$.  The
hypotheses imply that $\phi$ induces a linear surjection from the space
$\Lambda(K)^{1}$ of homogeneous degree $1$ elements of $\Lambda(K)$ to the
space $\Lambda \left( \partial \Delta^{n} \right)^{1}$ of homogeneous
degree $1$ elements of $\Lambda\left( \partial \Delta^{n} \right)$, which
is $(n+1)-$dimensional.  Thus there is a collection of $n+1$ elements of
the standard basis $\{ v_1, \ldots, v_m \}$ for $\Lambda(K)$ which map
onto a basis for the space of homogeneous degree $1$ elements of $\Lambda
\left( \partial \Delta^{n} \right)$.  We can thus assume, without loss of
generality, that $\left\{ \phi\left(v_1 \right) , \ldots, \phi\left(
v_{n+1} \right) \right\}$ is a basis for $\Lambda \left( \partial
\Delta^{n} \right)^{1}$.

For $i \in \{ n+2, \ldots, m \}$, let $s_i$ denote the (unique) linear
combination of $v_1, \ldots, v_{n+1}$ such that $\phi \left( s_i \right) =
\phi \left( v_i \right)$.  Thus, each $s_i + v_i$ is an element of
$\ker\phi$.  Since the set $\left\{ s_i + v_i \mid i \in \left\{ n+2,
\ldots, m \right\} \right\}$ is linearly independent, it must form a basis
for $\ker\phi \cap \Lambda(K)^{1}$, since the dimension of $\ker\phi \cap
\Lambda(K)^{1}$ is $m-n-1$.

Now assume that $\ker\phi$ is generated by degree $1$ and degree $2$
elements.  Suppose that $\ker\phi \cap \Lambda(K)^{2}$ is spanned by $t_1,
t_2, \ldots, t_k$, where $t_i$, for $i \in \left\{ 1, \ldots, k \right\}$,
is a homogeneous element of degree $2$.  
By Lemma \ref{flaglemma}, 
$$ \Lambda(K) \cong 
\Lambda \left[ v_1 , \ldots, v_m \right]/ \langle u_1, \ldots, u_l \rangle,$$
where $u_i$ is a homogeneous element of degree $2$ for $1 \leq i \leq l$.  
It follows that
  $$ \phi: \Lambda \left[ v_1 , \ldots, v_m \right] / \left\langle s_{n+2} + v_{n+2}, \ldots, s_{m} + 
  v_{m}, t_1, \ldots t_k, u_1 , \ldots , u_l  \right\rangle \rightarrow \Lambda \left( \partial 
  \Delta^{n} \right)$$
is an isomorphism.  
Let $\Lambda_{\phi} (K)$ be the quotient of 
$\Lambda \left[ v_1 , \ldots , v_{n+1} \right]$
by the ideal $I_{\phi}(K) = \langle \hat{t}_{1}, \ldots, \hat{t}_{k},
\hat{u}_{1}, \ldots, \hat{u}_{l} \rangle$, 
where $\hat{t}_{j}$ (respectively, $\hat{u}_{j}$)
is the result of replacing $v_i$ with $s_i$
($n+2 \leq i \leq m$) in $t_{j}$ (respectively, $u_j$). 
Note that $I_{\phi}(K)$ is generated by
homogeneous elements of degree $2$.  It is easy to see that the map
  $$\psi : \Lambda_{\phi}(K) \rightarrow \Lambda \left[ v_1 , \ldots , v_m \right]   
  / \left\langle s_{n+2} 
  + v_{n+2}, \ldots, s_{m} + v_{m}, t_1, \ldots t_k, u_1 , \ldots, u_l
  \right\rangle,$$
sending $v_i$ to $v_i$ for $i \in \{ 1, \ldots, n+1 \}$, is an 
isomorphism, which also preserves degree.

Now we obtain a contradiction by counting the dimensions of
$\Lambda_{\phi}(K)^{2}$, $\Lambda_{\phi}(K)^{n+1}$, $\Lambda\left(
\partial \Delta^{n} \right)^{2}$, and $\Lambda\left( \partial \Delta^{n}
\right)^{n+1}$ as vector spaces.  We have:
  $$\mathrm{dim} \left( \Lambda\left( \partial \Delta^{n} \right)^{2} 
  \right) = \frac{n(n+1)}{2};  \quad \quad \mathrm{dim} \left( 
  \Lambda\left( \partial \Delta^{n} \right)^{n+1} \right) = 0.$$ 
Either $I_{\phi}(K)$ is the $0$ ideal or it isn't.  If it is, then
  $$\mathrm{dim} \left( \Lambda_{\phi}(K)^{n+1} \right) = 1;$$
if it isn't, then
  $$\mathrm{dim} \left( \Lambda_{\phi}(K)^{2} \right) <  
  \frac{n(n+1)}{2}.$$

In either case, we have a contradiction since $\phi \circ \psi$ is a
degree-preserving bijection and thus preserves the dimension in each
degree.
\end{proof}

\begin{corollary} \label{cor:principle2}
Let $K_1$ and $K_2$ be finite simplicial complexes.  
\begin{enumerate}
\item If $\phi: \Lambda(K_1) \rightarrow \Lambda(K_2)$ is a
degree-preserving surjection, $K_1$ is a flag complex, and $\ker\phi$ is
generated by homogeneous elements of degrees one and two, then $K_2$ is
also a flag complex.

\item If $\phi: \Lambda(K_1) \rightarrow \Lambda(K_2)$ is a
degree-preserving isomorphism, then $K_1$ is a flag complex if and only if
$K_2$ is.
\end{enumerate}
\end{corollary}
\begin{proof}
(1) If $K_2$ is not flag, then for some $n \geq 2$, $\partial \Delta^{n}$
is a full subcomplex of $K_2$, i.e., $\partial \Delta^{n}$ is not the
boundary of an $n$-simplex in $K_2$.  Define a map $\psi: \Lambda(K_2)
\rightarrow \Lambda(\partial \Delta^{n})$, sending a given vertex $v$ to
$0$ if $v \not \in \partial \Delta^{n}$, and to itself otherwise.  The map
$\psi$ is a degree-preserving surjection whose kernel is generated by
elements of degree $1$.  It follows that $\psi \circ \phi: \Lambda(K_1)
\rightarrow \Lambda(\partial \Delta^{n})$ is a degree-preserving
surjection whose kernel is generated by homogeneous elements of degrees
$1$ and $2$.  This contradicts Proposition \ref{prop:principle}.

(2) This is an easy consequence of (1). 
\end{proof}

Note that Gubeladze \cite{Gubeladze} has proven a strong generalization 
of Corollary \ref{cor:principle2}(2): $\phi: \Lambda(K_1) \rightarrow 
\Lambda(K_2)$ is a degree-preserving isomorphism if and only if $K_1$ is 
isomorphic to $K_2$ as simplicial complexes.

Our experience in computing $H^{\ast}\left( \ud{n}{T}; \mathbb{Z}/
2\mathbb{Z} \right)$ for various small examples, including the case $T =
T_{min}$ and $n=4$, suggests the following conjecture:

\begin{conjecture}
The cohomology ring $H^{\ast} \left( \ud{n}{T} ; \mathbb{Z} / 2\mathbb{Z}
\right)$ is an exterior face algebra, for any tree $T$ and any $n$.
\end{conjecture}

We note finally that the conjecture seems just as likely to be true for
arbitrary fields, not simply $\mathbb{Z} / 2\mathbb{Z}$.

\subsection{Which Tree Braid Groups Are Right-Angled Artin?}

Recall the definition of a right-angled Artin group from Example
\ref{ex:RAAG}.  In this subsection, we characterize exactly which tree
braid groups are right-angled Artin.  Theorem \ref{thm:Ghrist} states that
a tree braid group $B_n T$ is a right-angled Artin group exactly when
either $n < 4$ or $T$ is linear (recall that a tree is \emph{linear} if
there exists an embedded line segment which contains all of the 
essential vertices of $T$).

Let $T_{min}$ be the minimal nonlinear tree described in Subsection
\ref{sec:computation}.

\begin{lemma} \label{lem:embed}
Let $n \geq 4$.  Let $T$ be a nonlinear tree that is sufficiently
subdivided for $n$.  There is a cellular embedding $\theta$ of (a suitably
subdivided)  $T_{min}$ into $T$ such that:
\begin{enumerate}
\item the image of $T_{min}$ is sufficiently subdivided for $4$ strands;  

\item there is a choice of basepoints $\overline{\ast}$ for $T_{min}$ 
and $\ast$ for $T$ such that $\overline{\ast}$ has degree $1$ in 
$T_{min}$, $\ast$ has degree $1$ in $T$, and the geodesic segment 
$\left[ \ast, \overline{\ast} \right]$ in $T$ crosses exactly $n-4$ 
edges, none of which are edges of $T_{min}$.
\end{enumerate}
\end{lemma}

\begin{proof}
Choose a collection $\mathcal{C}$ of essential vertices of $T$ such that
the elements of $\mathcal{C}$ all lie along an embedded arc, and such that
$\mathcal{C}$ is a maximal set of essential vertices with this property.  
Fix an arc $\left[ v_1, v_2 \right]$ in $T$ such that $\mathcal{C}
\subseteq \left[ v_1, v_2 \right]$, where $v_1$ and $v_2$ are essential.  
Since $T$ is a nonlinear tree, there exists an essential vertex $v_3 \not
\in \left[ v_1, v_2 \right]$. Consider the geodesic segment $\gamma$
connecting $v_3$ to $\left[ v_1, v_2 \right]$.  By maximality of
$\mathcal{C}$, $\gamma$ must meet $\left[ v_1, v_2 \right]$ in another
essential vertex $v_4 \not \in \left\{ v_1, v_2 \right\}$, for otherwise
$\gamma \cup \left[ v_1, v_2 \right]$ is an arc containing $\mathcal{C}
\cup \left\{ v_3 \right\}$, which contradicts the maximality of
$\mathcal{C}$.

\begin{figure}[!h]
\centering
\input{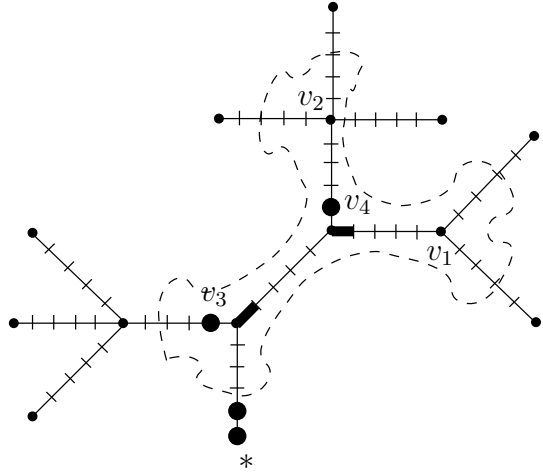}

\caption{The larger tree $T$ is sufficiently subdivided for $n=6$. The
smaller tree (encircled) is a copy of $T_{min}$.  This figure shows the
image of a critical $2$-cell in $\ud{4}{T_{min}}$ under the map $\theta :
\ud{4}{T_{min}} \rightarrow \ud{6}{T}$.}
\label{fig:embed}
\end{figure}

The $Y$-graph formed by the segments $\left[ v_1, v_4 \right]$, $\left[ 
v_2, v_4 \right]$, and $\left[ v_3, v_4 \right]$ is sufficiently 
subdivided for $4$, since the tree $T$ is sufficiently subdivided for 
$n$ and $n \geq 4$.  For $i=1,2,3$, add to the $Y$-graph two additional 
embedded line segments at $v_i$, each consisting of exactly $3$ edges, 
in such a way that the new segments have no edges in common with either 
each other or with the $Y$-graph.  It is possible to do this because 
each of the vertices $v_1$, $v_2$, and $v_3$ are essential.  The result 
of this procedure gives a cellular embedding of $T_{min}$ into $T$ which 
satisfies (1).

To produce an embedding satisfying (2) as well, proceed as follows.  
Choose a vertex $\hat{\ast}$ having degree $1$ in $T_{min}$.  If 
$\hat{\ast}$ has degree $1$ in $T$, then since $T$ is sufficiently 
subdivided for $n$ it must be that $n = 4$.  In this case, choose $\ast 
= \overline{\ast} = \hat{\ast}$.  Otherwise, $\hat{\ast}$ has degree at 
least $2$ in $T$.  Choose an arc $\hat{\gamma}$ in $T$ with no edges in 
common with the embedding of $T_{min}$, and with the embedding of 
$\hat{\ast}$ as one of its endpoints.  Furthermore, choose 
$\hat{\gamma}$ to be a maximal such arc, so that the other endpoint of 
$\hat{\gamma}$ has degree $1$ in $T$.  Declare the other endpoint of 
$\hat{\gamma}$ to be $\ast$, and let $\overline{\ast}$ be the (unique) 
vertex lying on $\hat{\gamma}$ at distance exactly $n-4$ from $\ast$.  
Modify $T_{min}$ (if necessary) by adding in the segment $\left[ 
\overline{\ast}, \hat{\ast} \right]$.
\end{proof}            

Let $n \geq 4$, and let $T$ be a nonlinear tree.  The embedding $\theta$
of the previous lemma induces a map of configuration spaces $\theta:
\ud{4}{T_{min}} \rightarrow \ud{n}{T}$, defined by $\theta \left(
\left\{ c_1, c_2, c_3, c_4 \right\} \right) = \left\{ \theta(c_1),
\theta(c_2), \theta(c_3), \theta(c_4) \right\} \cup \left\{ \ast, v_1,
v_2, \ldots, v_{n-5} \right\}$, where $\ast, v_1, v_2, \ldots, v_{n-5}$
are the $n-4$ vertices of $T$ closest to $\ast$ (see Figure
\ref{fig:embed}).  

Choose an embedding of $T$ in the plane, and consider the induced 
classifications of cells in $\ud{n}{T}$ and $\ud{4}{T_{min}}$ into 
the critical, collapsible, and redundant cell types.

\begin{proposition}
The map $\theta$ preserves cell type - i.e. takes critical, 
collapsible, and redundant cells in $\ud{4}{T_{min}}$ to critical, 
collapsible, and redundant cells in $\ud{n}{T}$, respectively.
\end{proposition}

\begin{proof}
Let $c$ be a cell in $\ud{4}{T_{min}}$.  Since the map $\theta$ on 
configuration spaces is induced by a cellular embedding on the level of 
trees, by the choices of $\ast$ and $\overline{\ast}$, a vertex in $c$ 
is blocked if and only if it is blocked in $\theta(c)$.  Since the 
embedding of $T_{min}$ in the plane is induced from the embedding of $T$ 
in the plane, an edge in $c$ is respectful if and only if it is 
respectful in $\theta(c)$.  The cell $\theta(c)$ is obtained from 
$c$ by adding exactly $n-4$ blocked vertices.  Thus, the numbering on 
vertices and edges in $\theta(c)$ used to determine cell type (see the 
discussion preceding Definition \ref{def:critical}) differs from the 
numbering for $c$ only by the insertion of $n-4$ blocked vertices at the 
beginning of the numbering.  By the definition of the Morse matching, 
Definition \ref{def:critical}, the proposition is proven.
\end{proof}

\begin{proposition} \label{prop:cohom} We have:
\begin{enumerate}
\item The map $\theta$ induces an injection $\theta_{\ast}: H_{\ast}
\left( \ud{4}{T_{min}} ; \mathbb{Z} / 2\mathbb{Z} \right)  \rightarrow
H_{\ast} \left( \ud{n}{T} ; \mathbb{Z} / 2\mathbb{Z} \right)$.  The
homology class corresponding to a given critical cell $c \subseteq
\ud{4}{T_{min}}$ goes to a homology class corresponding to $\theta(c)
\subseteq \ud{n}{T}$.  In particular, the image of $\theta_{\ast}$ is a
direct factor of $H_{\ast}\left( \ud{n}{T} \right)$.

\item The induced map $\theta^{\ast}: H^{\ast}(\ud{n}{T};
\mathbb{Z}/2\mathbb{Z}) \rightarrow H^{\ast}(\ud{4}{T_{min}};
\mathbb{Z}/2\mathbb{Z})$ sends the dual of a critical cell $c$ to
$(\theta^{-1}(c))^{\ast}$ if $c$ is in the image of $\theta$, or to $0$
otherwise.
\end{enumerate}
\end{proposition}

\begin{proof} 
(1)  Fix an embedding of $T$ into the plane, and choose an embedding
$\theta: T_{min} \rightarrow T$ as in Lemma \ref{lem:embed}.  We note
that, due to the choices of the embedding $\theta: T_{min} \rightarrow T$
and basepoints, the map $\theta: \ud{4}{T_{min}} \rightarrow \ud{n}{T}$
sends collapsible cells to collapsible cells, redundant cells to redundant
cells, and critical cells to critical cells.

If $c$ is an arbitrary critical cell of $\ud{4}{T_{min}}$, then a cycle
representing the homology class determined by $c$ is $f^{\infty}(c)$,
which has the form $c + (collapsible~cells)$.  Since $\theta$ preserves a
cell's type, it follows that the homology class $\theta_{\ast}(c)$ may be
represented by a cycle of the form $\theta(c) + (collapsible~cells)$,
where $\theta(c)$ is critical.  By Lemma \ref{lem:homology}(2), the cycle
$\theta(c) + (collapsible~cells)$ is homologous to $f^{\infty}(
\theta(c) + (collapsible~cells))$.  By Lemma \ref{lem:homology}(3) and
the fact that $\theta(c)$ is critical,
  $$f^{\infty}(\theta(c) + collapsible~cells) = \theta(c) + 
  ((different)~collapsible~cells).$$

On the other hand, the homology class corresponding to the critical cell
$\theta(c)$ is, by definition, $f^{\infty}(\theta(c))$, which consists of
$\theta(c) + (collapsible~cells)$.  Thus, using the fact from Lemma
\ref{lem:homology}(3) that an $f$-invariant chain is determined by its
critical cells, we conclude that $\theta_{\ast}(c) = \theta(c)$, as
required.

(2) This is an easy consequence of (1) and the naturality of the universal
coefficient isomorphism.
\end{proof}

\begin{theorem} \label{thm:Ghrist}
The tree braid group $\bn{T}$ is a right-angled Artin group if and only
if $T$ is linear or $n < 4$.
\end{theorem}

\begin{proof}
($\Leftarrow$)  Connolly and Doig \cite{ConnollyDoig} showed that $B_{n}T$
is a right-angled Artin group if $T$ is linear.  If $n<4$ and $T$ is a
tree, then Theorem 4.3 of \cite{FS1} shows that $B_{n}T$ is in fact a free
group, since $\ud{n}{T}$ strong deformation retracts on a graph.  This
proves one direction.

(Note that it is possible to get a proof of Connolly and Doig's result as
an application of the ideas in \cite{FS1}.  Suppose $T$ is a linear tree.  
Choose some basepoint $\ast$ for $T$ and an embedded arc $\ell$ containing
$\ast$ and all essential vertices of $T$; it is possible to do this since
$T$ is linear.  Now embed $T$ in $\mathbb{R}^{2}$ so that: (1) $\ast$ is
mapped to the origin; (2) $\ell$ is mapped to a segment on the positive
$y$-axis, and (3) the image of $T$ is contained in $\{ (x,y) \in
\mathbb{R}^{2} \mid x \leq 0 \}$. With this choice of embedding and the
induced order on the vertices of $T$, Theorem 5.3 of \cite{FS1} gives a
presentation of $B_{n}T$ as a right-angled Artin group.  The proof is left
as an exercise for the interested reader.)

($\Rightarrow$)  Proof by contradiction.  Suppose that $T$ is
nonlinear, $n \geq 4$, and $B_{n}T$ is a right-angled Artin group.  Since
$\ud{n}{T}$ is aspherical ( \cite{Abrams}, \cite{Ghrist}), $\ud{n}{T}$ is
a $K( B_{n}T, 1)$. In particular, by Proposition \ref{prop:flag}, the
cohomology ring $H^{\ast}(\ud{n}{T}; \mathbb{Z}/2\mathbb{Z})$ is the
exterior face algebra of a flag complex.

We choose an embedding $\theta: \ud{4}{T_{min}} \rightarrow \ud{n}{T}$ as
in Lemma \ref{lem:embed}.  By Proposition \ref{prop:cohom},
$\theta^{\ast}: H^{\ast}(\ud{n}{T}; \mathbb{Z}/2\mathbb{Z}) \rightarrow
H^{\ast}(\ud{4}{T_{min}}; \mathbb{Z}/2\mathbb{Z})$ is surjective, and it
is necessarily degree-preserving.  Since $H^{\ast}(\ud{4}{T_{min}};
\mathbb{Z}/2\mathbb{Z})$ is the exterior face algebra of a complex that is
not flag, we will arrive at a contradiction to Corollary
\ref{cor:principle2}(1) if we can show that $\ker(\theta^{\ast})$ is
generated by homogeneous elements of degrees one and two.  For this, it is
sufficient to show that if $c$ is a critical cell in $\ud{n}{T}$ of
dimension at least $3$, then $c^{\ast}$ is divisible by some element
$c_{1}^{\ast} \in \ker(\theta^{\ast})$ of degree one.

Let $c$ be a critical cell in $\ud{n}{T}$ of dimension at least $3$.  
There are two cases: either every cell of $c$ lies inside of (the embedded
image of) $T_{min} \cup [\ast,\overline{\ast}]$, or some cell of $c$ is 
not contained in $T_{min} \cup [\ast,\overline{\ast}]$.

We first consider the case in which some vertex or edge $x$ of $T$
occurring in $c$ is not contained in $T_{min} \cup
[\ast,\overline{\ast}]$.  Either $x$ is an edge $e$ or $x$ is a blocked
vertex.  If $x$ is a blocked vertex, then at the largest essential vertex
on the geodesic $[x, \ast]$ there must be a disrespectful edge $e$.  
In either case, break all edges of $c$ other than $e$, and consider the
resulting $1$-cell $c'$.  By Lemma \ref{lem:big}(4), $c'$ is equivalent to
a critical $1$-cell $\tilde{c'}$, and the proof of Lemma \ref{lem:big}(4)
shows that $\tilde{c'}$ may be described as simply the result of moving
all vertices in $c'$ toward $\ast$ until they are all blocked.  If follows
that $x$ occurs in $\tilde{c'}$. This implies that $\tilde{c'}$ is not in
the image of $\theta: \ud{4}{T_{min}} \rightarrow \ud{n}{T}$, since all
cells in this image consist of cells in $T_{min}$.  It now follows from
Proposition \ref{prop:cohom}(2) that $\tilde{c'}^{\ast} \in
\ker(\theta^{\ast})$.  But $[c]$ is the least upper bound of its
$1$-dimensional lower bounds, so Proposition \ref{prop:cup}(2) implies
that $\tilde{c'}^{\ast} \mid c^{\ast}$, as required.

Finally, suppose that all vertices and edges in $c$ are contained in 
$T_{min} \cup [\ast, \overline{\ast}]$.  Let $A$, $B$, $C$, $D$ denote 
the four essential vertices of $T_{min}$, listed in the order they are 
numbered, from least to greatest.  Let $e_{A}$, $e_{B}$, $e_{C}$, and 
$e_{D}$ denote the edges incident with $A$, $B$, $C$, and $D$, 
respectively, which are in the greatest direction possible from each - 
namely, 2 (see the paragraphs preceding Definition \ref{def:critical}). 
Note these are the only four edges in $T_{min}$ which can possibly be 
disrespectful in a cell of $\ud{n}{T_{min}}$.  Since $c$ is 
critical and has dimension at least $3$ by assumption, $c$ contains at 
least $3$ edges in $T_{min}$, and these must be chosen from $\{ e_{A}, 
e_{B}, e_{C}, e_{D} \}$.  It follows that either $e_{A}$ or $e_{B}$ is 
in $c$.  Let $c'$ be the result replacing all edges in $c$ with either 
endpoint, except for $e_{A}$ if $e_{A} \in c$ or $e_{B}$ if $e_{A} \not 
\in c$.  Let $\tilde{c'}$ be the result of moving all vertices in $c'$ 
toward $\ast$ until they are blocked.  By Lemma \ref{lem:big}(4), 
$\tilde{c'}$ is critical.

For an arbitrary cell $\overline{c}$, a vertex $v \in \overline{c}$ is 
\emph{blocked by} an edge $e \in \overline{c}$ if and only if there are 
no vertices of $T - \overline{c}$ which are less than $v$ and between 
$v$ and $e$.  We claim that there are at least five vertices in 
$\tilde{c'}$ blocked by $e_A$ (if $e_A \in c$) or $e_B$ (if $e_A \not\in 
c$).  The reason is that there must be at least three edges in $c$, all 
of which are greater than or equal to $e_A$ or $e_B$, respectively.  As 
each edge must be disrespectful, each of the three edges blocks 
at least one vertex.  When all of these strands are moved towards $\ast$ 
until they are blocked, the result is that at least five vertices
are blocked in $\tilde{c'}$ by the edge $e_{A}$ 
or $e_{B}$, respectively.  This proves the claim.

It follows from this that $\tilde{c'}$ is not in the image of $\theta:
\ud{4}{T_{min}} \rightarrow \ud{n}{T}$, since any cell in
$\theta(\ud{4}{T_{min}})$ will contain at most $4$ cells from the tree
$T_{min}$.  This implies that $(\tilde{c'})^{\ast} \in
\ker(\theta^{\ast})$.  By construction $[\tilde{c'}] \leq [c]$, and, since
$[c]$ is the least upper bound of its $1$-dimensional lower bounds,
$(\tilde{c'})^{\ast} \mid c^{\ast}$ by Proposition \ref{prop:cup}(2).
\end{proof}


\bibliography{refs-FS2}
\nocite{*}
\bibliographystyle{plain}

\end{document}